\documentclass[pdflatex,sn-mathphys-num]{sn-jnl}

\usepackage{graphicx}%
\usepackage{multirow}%
\usepackage{amsmath,amssymb,amsfonts}%
\usepackage{amsthm}%
\usepackage{mathrsfs}%
\usepackage[title]{appendix}%
\usepackage{xcolor}%
\usepackage{textcomp}%
\usepackage{manyfoot}%
\usepackage{booktabs}%
\usepackage{algorithm}%
\usepackage{algorithmicx}%
\usepackage{algpseudocode}%
\usepackage{listings}%

\usepackage{amsmath}
\usepackage{amsthm}
\usepackage{amssymb}
\usepackage{graphicx}
\usepackage{standalone}
\usepackage{color}
\usepackage{comment}
\usepackage{enumerate}
\usepackage[normalem]{ulem}

\theoremstyle{plain}
\newtheorem{thm}{Theorem}[section]

\newtheorem{lem}[thm]{Lemma}
\newtheorem{prop}[thm]{Proposition}
\newtheorem{defn}[thm]{Definition}
\newtheorem{cor}[thm]{Corollary}
\theoremstyle{definition}
\newtheorem{rem}[thm]{Remark}

\newtheorem{construction}[thm]{Construction}
\newtheorem{problem}{Open Problem}

\title{Regular fat linear sets}
\date{\today}

\newcommand{\F}{\mathbb F}
\newcommand{\Fq}{{{\mathbb F}_q}}
\newcommand{\fq}{{{\mathbb F}_q}}

\newcommand{\Fqn}{{{\mathbb F}_{q^n}}}
\newcommand{\fqn}{{{\mathbb F}_{q^n}}}
\newcommand{\Fqt}{{{\mathbb F}_{q^t}}}

\newcommand{\la}{\langle}
\newcommand{\ra}{\rangle}

\newcommand{\C}{\mathcal C}
\newcommand{\Rt}{\operatorname{R-}q^t\operatorname-}

\newcommand{\eps}{\varepsilon}
\newcommand{\ii}{^\iota}

\DeclareMathOperator{\Mod}{mod}
\DeclareMathOperator{\PG}{PG}
\DeclareMathOperator{\N}{N}
\DeclareMathOperator{\GaL}{\Gamma L}
\DeclareMathOperator{\GL}{GL}
\DeclareMathOperator{\PGL}{PGL}

\DeclareMathOperator{\Tr}{Tr}

\begin{document}

\title[Regular fat linear sets]{Regular fat linear sets}


\author[1]{\fnm{Valentino} \sur{Smaldore}}\email{valentino.smaldore@unipd.it}

\author*[1]{\fnm{Corrado} \sur{Zanella}}\email{corrado.zanella@unipd.it}

\author[2]{\fnm{Ferdinando} \sur{Zullo}}\email{ferdinando.zullo@unicampania.it}

\affil[1]{\orgdiv{Dipartimento di Tecnica e Gestione dei Sistemi Industriali}, \orgname{Universit\`a degli Studi di Padova}, \orgaddress{\street{Stradella S. Nicola, 3}, \city{Vicenza VI}, \postcode{35030}, \country{Italy}}}

\affil[2]{\orgdiv{Dipartimento di Matematica e Fisica}, \orgname{Universit\`a degli Studi della Campania ``Luigi Vanvitelli''}, \orgaddress{\street{Viale Lincoln 5}, \city{Caserta CE}, \postcode{81100}, \country{Italy}}}

\abstract{In this work, we introduce $(r,i)$-regular fat linear sets, which are defined as linear sets containing exactly $r$ points of weight $i$ and all other points of weight one. This notion generalizes and unifies existing constructions; scattered linear sets, clubs, and other previously studied families are special cases. We present new classes of regular fat linear sets in $\PG(k-1,q^n)$ for composite $n$ and study their equivalence classes. Finally, we show that regular fat linear sets naturally yield three-weight rank-metric codes, which we use to obtain bounds on their parameters.}

\keywords{Linear set, Linearized polynomial, Rank-metric code, Scattered polynomial}

\pacs[MSC Classification]{51E20, 11T06, 11T71, 94B05}

\maketitle

\section{Introduction}
In~\cite{Lunardon}, Lunardon formally introduced linear sets as a tool for constructing a special type of blocking set in the projective plane. 
Linear sets are defined as sets of points \(L_U\) in the projective space \(\PG(k-1,q^n)\), where the points are determined by an $\Fq$-subspace \(U\) of \(\F_{q^n}^k\). 
Subgeometries are a basic class of examples, and many further constructions have been developed over the years.

Linear sets are highly useful in various classification results and constructions in finite geometry and coding theory due to their deep connections with algebraic and geometric structures (see~\cite{lavrauw2015field,polverino2010linear}). 
A central notion associated with linear sets is the \emph{weight} of a point, which intuitively measures how much the point overlaps the subspace \(U\).

In~\cite{BBL2000}, Ball, Blokhuis, and Lavrauw  constructed a special type of linear set in which all points have weight one. 
Later, in~\cite{BlLa00}, Blokhuis and Lavrauw initiated the study of these linear sets, calling them \emph{scattered linear sets}, and investigating them in a broader context. 
Recently, these sets have  attracted considerable interest due to their connections with maximum rank distance codes and, more generally, rank-metric codes. See~\cite{polverino2020connections} for a survey.

Fancsali and Sziklai in~\cite{FanSzi06} (see also~\cite{FanSzi09}) studied maximal partial $2$-spreads in $\PG(8,q)$, which are closely related to linear sets on the projective line in which all but one point have weight one. 
If the unique point of weight greater than one has weight $i$, then the set is called an $i$-\emph{club}. 
In particular, clubs in even characteristic are linked to translation KM-arcs, as shown by De Boeck and Van de Voorde in~\cite{DBoVdV16}. More recently, Alfarano, Byrne, and Fulcher in \cite{ABF24} showed that the existence of clubs is related to the representability of the free product of two $q$-uniform matroids. See \cite{NapPolSanZul24} for other applications.
A recent investigation of clubs in projective spaces can be found in \cite{MSZ25}.
Additionally, \cite{NapPolSanZul22} investigates linear sets with precisely two points of weight greater than one. This topic is later expanded upon in~\cite{Zul23,santonastaso}. 
The existence of these linear sets relates to the representability problem of the direct sum of uniform $q$-matroids, as proven in~\cite{alfarano}.

In this paper, we introduce a natural generalization. We define an $(r,i)$-\emph{regular fat  linear set} as a linear set with $r$ points of weight $i$  and all other points of weight one. An $(r,i)$-\emph{regular fat linearized polynomial} is one for which the related linear set of rank $n$ in $\PG(1,q^n)$ is $(r,i)$-regular fat. This notion extends previous definitions: a $(0,i)$-regular fat linear set is a scattered linear set; a $(1,i)$-regular fat linear set is an $i$-club; and a $(2,i)$-regular fat linear set of rank $2i$ corresponds to some of the linear sets studied in \cite{NapPolSanZul22}.
In \cite{BMZZ}, the authors defined an $r$-fat linearized polynomial as one whose associated linear set has exactly $r$ points of weight greater than one. 
It seems that, when $r$ is small, finding examples is very difficult. In fact, the authors of  \cite{BMZZ} proved that there are no exceptional $r$-fat polynomials, i.e., linearized polynomials that remain $r$-fat over infinitely many extensions of the field coordinatizing the projective line. In this paper, we 
review the known examples of regular fat linear sets and polynomials and derive some properties and introduce two classes of $(r,t)$-regular fat linear sets in $\mathrm{PG}(k-1,q^n)$ when $n$ is a multiple of $t$. Consequently, we examine the equivalence problem among linear sets belonging to these classes.
Next, we study the polynomial
\begin{equation}\label{q:pol}
\phi_{m,\sigma}=X^{{\sigma}^{t-1}}+X^{{\sigma}^{2t-1}}+m(X^{\sigma}-X^{{\sigma}^{t+1}})\in\mathbb{F}_{q^{2t}}[X],
\end{equation}
introduced in \cite{SmZaZu2},
where $q$ is odd, $\sigma$ is a generator of $\mathrm{Gal}(\fqn\colon \fq)$, and $t\ge3$. We characterize its behavior as a regular fat polynomial.
Specifically, when $m$ satisfies certain conditions, the associated linear set belongs to the aforementioned class of regular fat linear sets.

Following Randrianarisoa's paper \cite{Randrianarisoa2020ageometric}, the relationship between linear sets (more precisely, $q$-systems) and rank-metric codes was examined in depth, yielding insights into both fields.
In fact, the optimality of rank-metric codes has been fully characterized from a geometric perspective.
As with the Hamming metric, recent attention has been given to  rank-metric codes with few weights. A complete classification of one-weight codes was provided in \cite{Randrianarisoa2020ageometric} (see also \cite{lincutt}), and a geometric description of two-weight rank-metric codes was given in \cite{ZulPolSanShe24} via the scatteredness property. A complete geometric characterization of the three-weight case is still missing, though several examples have been presented in \cite{MSZ25}.
We demonstrate that, via a duality argument, regular fat linear sets provide examples of three-weight rank-metric codes. This extends the classes introduced in \cite{MSZ25}. Furthermore, we use the MacWilliams identities for rank-metric codes to derive a nontrivial bound on the number of points of weight greater than one in a regular fat linear set, leveraging the aforementioned connection between regular fat linear sets and three-weight rank-metric codes.
We conclude the paper by outlining possible future research directions that may be of interest.

\medskip
\subsection{Organization of the paper}

This paper is structured as follows: 
Section~\ref{s:regular} defines regular fat linear sets and reviews those present in the literature. 
To our knowledge, no $(r,i)$-regular fat linear set with $r,i>2$ is known.
Motivated by this, in Section~\ref{s:examples}, we present two constructions of $((q^k-1)/(q-1), t)$-regular fat linear sets in the projective space $\PG(k-1, q^{\ell t})$, where $\ell \mid q^{t}-1$.
We deduce a polynomial form (see Theorem~\ref{t:polform}) when $k=2$ and the rank is $\ell t$. We also describe the equivalence classes of those sets under the action of $\GaL(k,q^{\ell t})$. Section~\ref{s:4} proves that the polynomial \eqref{q:pol} defined in \cite{SmZaZu2} is regular fat for any $m\in\F_{q^t}^*$ (see Theorem~\ref{t:summary}).  Section~\ref{s:5} associates a rank-metric code with each regular fat linear set. We were able to prove relationships between relevant parameters using some results from the theory of rank-metric codes. Section~\ref{s:conclusions} contains remarks on parameter $r$ and poses two open problems.

\subsection{Notation}

Let $n,k>1$ be integers and $q$ be a prime power.  
We use the notation $(\alpha_1:\alpha_2:\ldots:\alpha_k)$ to represent the point $\la(\alpha_1,\alpha_2,\ldots,\alpha_k)\ra_{\Fqn}$ in $\PG(k-1,q^n)$.
The \emph{$\Fq$-linear set} (or simply \emph{linear set}) \emph{of rank $\rho$} in $\PG(k-1,q^n)$ associated with an $\Fq$-subspace $U$ of $\F_{q^n}^k$, $\dim_{\Fq}U=\rho$, is
        \[
    L_U=\{\la v\ra_{\Fqn}\colon v\in U, v\neq 0\}.
        \]
The \emph{weight} with respect to $L_U$ of a point $P=\la v\ra_{\Fqn}\in\PG(k-1,q^n)$ is
\[
  w_{L_U}(P)=\dim_{\F_q}\left(\la v\ra_{\Fqn}\cap U\right).
\]
An $\Fq$-linear set
$L_U$ is \emph{scattered} \cite{BlLa00} if $\dim_{\F_q}\left(\la  v\ra_{{\Fqn}}\cap U\right)\le1$ for all $  v\in \F_{q^n}^k$.
For $1<t\mid n$, 
$L_U$ is \emph{$\Rt$partially scattered} if
$\dim_{\F_q}\left(\la  v\ra_{{\Fqt}}\cap U\right)\le1$ for all $  v\in \F_{q^n}^k$   \cite{LZ2}.

Consider the following set of $\fq$-\emph{linearized polynomials} (or $q$-\emph{polynomials}): 
\[
\mathcal L_{n,q}=\left\{\sum_{i=0}^{n-1}a_iX^{q^i}\colon a_0,\ldots,a_{n-1}\in\Fqn\right\}.
\]
As is well known, $\mathcal L_{n,q}$ is in one-to-one correspondence with the set of $\Fq$-linear endomorphisms of the $\Fq$-vector space $\Fqn$.
For any $f$ in $\mathcal L_{n,q}$ define $U_f=\{(x,f(x))\colon x\in\Fqn\}$. The \emph{linear set} of rank $n$ associated with $f$ is \[
  L_f=L_{U_f}=\{(x:f(x))\colon x \in\F_{q^n}^*\}\subseteq\PG(1,q^n).
\]
Conversely, every $\Fq$-linear set of rank $n$ in $\PG(1,q^n)$ can be represented, up to projectivities, in the form $L_f$, where $f$ is a linearized polynomial.

A linearized polynomial $f\in\mathcal L_{n,q}$ is called \emph{scattered} if $L_f$ is scattered. Equivalently,
\[ x,y\in\F_{q^n}^*,\quad
\frac{f(x)}x=\frac{f(y)}y\ \Rightarrow\ \frac xy\in\Fq.
\]      
Similarly, a linearized polynomial $f$ is called \emph{$\Rt$partially scattered} if $L_f$ is $\Rt$partially scattered.

In this paper, we use the following notation for the norm and trace of an element $x\in\Fqn$ over $\Fq$:
\begin{eqnarray*}
    \N_{q^n/q}(x)&=&x^{(q^n-1)/(q-1)},\\
    \Tr_{q^n/q}(x)&=&x^{q^{n-1}}+x^{q^{n-2}}+\cdots+x^q+x.
\end{eqnarray*}

Let $t\ge3$ be an integer.
Assume $\sigma=q^J$, $J\in\{1,\ldots,2t-1\}$ and $\gcd(J,2t)=1$.
According to \cite[Theorem 2.2]{SmZaZu2},  for any  $m\in\mathbb F_{q^t}^*$, the polynomial $ \phi_{m,\sigma}$ in \eqref{q:pol}
is $\Rt$partially scattered.
If, in addition, $m$ is neither a {$(\sigma-1)$}-th nor a $(\sigma+1)$-th power of an element of $\ker\operatorname{Tr}_{q^{2t}/q^t}$, then $\phi_{m,\sigma}$ is scattered.

\section{Regular fat linear sets and polynomials}\label{s:regular}

The precise definition of a regular fat linear set requires an additional condition on points of weight one, in order to ensure that it is indeed an $\Fq$-linear set. See also \cite[Theorem 2]{CsMaPe24} in this regard.

\begin{defn}
An \emph{$(r, i)$-regular fat  linear set}  ($r\ge0$, $i\ge2$) in $\PG(k-1,q^n)$ admits exactly $r$ points of weight $i$  and the remaining points have weight one, of which there exists at least one such point.

 An \emph{$(r, i)$-regular fat  linearized polynomial} is an $\fq$-linearized polynomial $f$ such that $L_f$ is an $(r,i)$-regular fat linear set {in $\PG(1,q^n)$}.
\end{defn}

The $(r,i)$-regular fat linear sets in the projective line are particular linear sets associated with $r$-fat linearized polynomials as defined in \cite{BMZZ}, where no requirement is set for points of weight greater than one.

\begin{rem} {For any integer $i\ge2$, a $(0,i)$-regular fat linear set is a scattered linear set, and vice versa.} \end{rem}

The size of a regular fat linear set is determined by its parameters. This can be proven using elementary equations that involve the number of points with a given weight, as described in \cite[(1-2)]{polverino2010linear}.
\begin{prop}
If $L_U$ is an $(r, i)$-regular fat  linear set of rank $\rho$, then 
\begin{equation}\label{e:size} |L_U|=\frac{q^\rho-1-r(q^i-q)}{q-1}.
\end{equation}
In particular, if $f$ is an $(r,i)$-regular fat $q$-polynomial, then
\[
|L_f|=q^{n-1}+q^{n-2}+\cdots+q^i+1-(r-1)(q^{i-1}+q^{i-2}+\cdots+q).
\]
\end{prop}

\subsection{Known examples for \texorpdfstring{$r=1$ or $i=2$}{r=1 or i=2}}
{For the remainder of this section, we will focus on the projective line. We review the literature on regular fat linear sets in $\PG(1,q^n)$.
The $(1,i)$-regular fat linear sets are known as \emph{$i$-clubs} in the literature and have been widely studied since the introduction of 2-clubs of rank 3 in $\PG(1,q^3)$  in \cite{FanSzi06,FanSzi09}.}
Their investigation was continued in \cite{LavVdV10}.
In \cite{LavZan15,LavZan16} the case $i=\rho-1$ is characterized and described in the Bruck-Bose setting.
In \cite{CsaMarPol18} it is proved that for $\rho=n$, $i=n-1$,  they only arise from the trace map.
In \cite{DBoVdV16} the authors deal with general $i$ in $\PG(1,q^n)$ and construct many examples.
Further results and constructions with $i=n-2$ can be found in \cite{NapPolSanZul24}.
More precisely:

\begin{itemize}
    \item The linear set $L_{\mathrm{Tr}_{q^n/q}}$ is an example of an $(n-1)$-club. In \cite[Theorem 3.7]{CsaMarPol18}, it was proven that every $(n-1)$-club is $\mathrm{PGL}(2,q^n)$-equivalent to $L_{\mathrm{Tr}_{q^n/q}}$.
    \item Let $n=\rho=\ell t$, $i=t(\ell-1)$, $\gcd(s,t)=1$ and $\sigma\colon x \in \fqn \mapsto x^{q^s} \in \fqn$. Then the linear set $L_T$, where 
\begin{equation}\label{eq:firstKMarc}
 T(x)=\left(\mathrm{Tr}_{q^{\ell t}/q^{t}}\circ\sigma\right)(x), 
\end{equation}
is an $i$-club in $\mathrm{PG}(1,q^{n})$ (see \cite[Theorem 3.3]{DBoVdV16}).
\item These are the only known examples of polynomials defining $i$-clubs, except for \cite[Corollary 5.5]{BMZZ} for $n=4$ (see also \cite{CsZ2018}) and \cite[Corollary 6.3]{BMZZ} for $n=5$.
\end{itemize}

In \cite{DBoVdV16,GW} two other examples of clubs are introduced.

\begin{construction}\label{constr:(n-2)clubpolynomial}
In \cite[Lemma 2.12]{DBoVdV16} it is proved that for any $\lambda \in \F_{q^n}^*$ such that $\{1,\lambda,\ldots,\lambda^{n-1}\}$ is an $\fq$-basis of $\fqn$, the $\fq$-linear set
\begin{equation}\label{eq:exKMnopol} L_{\lambda}=\{ (t_1\lambda+\cdots+t_{n-1}\lambda^{n-1}:t_{n-1}+t_n \lambda) \colon (t_1,\ldots,t_n)\in \mathbb{F}_q^n\setminus\{ {0}\} \},
\end{equation}
is an $(n-2)$-club of rank $n$ in $\PG(1,q^n)$.
\end{construction}

\begin{construction}\label{constr:clubviascattered}
In \cite[Lemma 3.6]{DBoVdV16} and \cite{GW} the following $i$-clubs are introduced.
Let $n=\ell t$, $t,\ell>1$ and $f \in \mathbb{F}_{q^t}[X]$ a scattered $q$-polynomial. Let 
\begin{equation}\label{eq:Ua,b}U_{a,b}=\left\{ \left(f(x_0)-ax_0,bx_0+\sum_{j=1}^{\ell-1}x_j \omega^j\right) \colon (x_0,x_1,\ldots,x_{\ell-1}) \in \F_{q^t}^\ell \right\},\end{equation}
for some fixed $a,b \in \F_{q^t}$ with $b\ne 0$ and $\{1,\omega,\ldots,\omega^{\ell-1}\}$ an $\F_{q^t}$-basis of $\F_{q^n}$. Then $L_{U_{a,b}}$ is an $i$-club of rank $n$,  with
\[ i=\left\{ \begin{array}{ll} t(\ell-1), & \text{if}\,\, f-aX\,\, \text{is invertible over}\,\,\F_{q^t},\\
t(\ell-1)+1, & \text{otherwise}.\end{array} \right. \]
\end{construction}

Consider the $q$-polynomial $g=X^{q^{s(n-1)}}+\delta X^{q^s} \in \fqn[X]$, where $s$ satisfies $\gcd(s,n)=1$ and $q$ is any prime power. The polynomial $g$ is called \emph{LP-polynomial} after Lunardon and Polverino who first introduced the $\fq$-linear set $L_{g}$ of $\PG(1,q^n)$.
Note that the polynomial $g$ satisfies
\begin{equation}\label{eq:fg} 
\dim_{\fq}(\ker(g-mX))\leq 2, 
\end{equation}
for all $m \in \fqn$ (see, for example, Theorem~\ref{t:GQ}).
Therefore, the points of $L_{g}$ have weight at most two, thus the number of points of weight two in $L_g$ is given by the following result.

\begin{thm}\label{th:finalLP}\cite[Theorem 7.12]{BMZZ}
Let $\delta\in\mathbb{F}_{q^n}$, $\gcd(s,n)=1$ and $g=X^{q^{s(n-1)}}+\delta X^{q^s} \in \fqn[X]$.
Then $L_g$ is an $(r,2)$-regular fat linear set where 
\begin{itemize}
    \item $r=0$, if $\N_{q^n/q}(\delta)\ne 1$;
    \item $r=\frac{q^{n-1}-1}{q^2-1}$, if $\N_{q^n/q}(\delta)= 1$ and $n$ is odd;
    \item $r=\frac{q^2(q^{\frac{n-2}2}+1)(q^{\frac{n-2}2-1}-1)}{q^2-1}+1$, if $\N_{q^n/q^2}(\delta)= 1$ and $n\equiv 0 \pmod{4}$;
    \item $r=\frac{q^2(q^{\frac{n-2}2-1}+1)(q^{\frac{n-2}2}-1)}{q^2-1}+1$, if $\N_{q^n/q^2}(\delta)= -1$ and $n\equiv 2 \pmod{4}$;
    \item $r=\frac{(q^{n/2-1}+1)(q^{n/2}-1)}{q^2-1}$, if $\N_{q^n/q}(\delta)= 1$, $\N_{q^n/q^2}(\delta)\ne  1$ and $n\equiv 0 \pmod{4}$;
    \item $r=\frac{(q^{n/2}+1)(q^{n/2-1}-1)}{q^2-1}$, if $\N_{q^n/q}(\delta)= 1$, $\N_{q^n/q^2}(\delta)\ne  -1$ and $n\equiv 2 \pmod{4}$.
\end{itemize}
\end{thm}

Another example of regular fat linearized polynomial is given by $f=X^q+\delta X^{q^2}$. Indeed, by \cite[Proposition 3.8 and Remark 3.9]{Zan}, we have that $L_f$ is $(r,2)$-regular fat with $r>0$.

Starting from a result in \cite{LavVdV10}, the authors in \cite{DBoVdV22} 
 find the weight distributions of the linear sets of rank $\rho\le4$ in $\PG(1,q^n)$ and rank 5 in $\PG(1,q^5)$. For $\rho=4$, they prove that all $(r,i)$-regular fat linear sets have the following parameters:
\[   (1,3),\quad (q+1,2),\quad (1,2), \quad \text{and }(2,2),\]
and all possibilities occur.

In \cite{CsaMarPolZ18} it is proved that all elements in a class of linear sets of type $L_{bX+X^{q^4}}\subseteq\PG(1,q^6)$ are either scattered or $(r,2)$-regular fat; in \cite{BarCsaMon21} the elements of both subclasses are characterized in terms of $b$.

Finally, the main result in \cite{PolZul20} leads to the construction of a large class of $(r,2)$-regular fat linearized polynomials.

\subsection{Known examples for \texorpdfstring{$r=2$}{r=2}}

The state of the art for $(2,i)$-regular fat linear set can be found in \cite{NapPolSanZul22}, 
where the notion of a \emph{linear set with complementary weights} is introduced. 
This is one that has two points such that the sum of the weights of the points equals the rank of the linear set.


\begin{thm}\cite[Proposition 3.2]{NapPolSanZul22}
    Let $L_W$ be an $\Fq$-linear set of rank $\rho$ in $\PG(1, q^n )$ for which there exist two distinct points $P , Q \in L_W$ such that
$w_{L_W} ( P ) = s$, $w_{L_W} ( Q ) = s'$ and $s + s' = \rho \le n$. Then, $L_W$ is $\PGL(2, q^n )$-equivalent to a linear set $L_U$  where $U = S \times S'$ for
some $\Fq$-subspaces $S$ and $S'$ of $\Fqn$, with $\dim_{\fq}( S ) = s$, $\dim_{\fq}( S' ) = s'$.
In addition, $S\cap S'=\{0\}$ can be assumed.
\end{thm}

\begin{thm}\cite[Theorem 3.4]{NapPolSanZul22}
    Let $L_W$ be an $\Fq$-linear set of rank $n$ in $\PG(1, q^n )$ for which there exist two distinct points $P , Q \in L_W$ such that
$w_{L_W} ( P ) = s$, $w_{L_W} ( Q ) = s'$ and $s + s' = n$. Then,  for
some $\Fq$-subspaces $S$ and $S'$ of $\Fqn$ with $\dim_{\fq}( S ) = s$, $\dim_\Fq( S' ) = s'$,
$\Fqn=S\oplus S'$, up to projectivities
\[
L_W=L_{p_{S,S'}}=\{(x:p_{S,S'}(x))\colon x\in\F_{q^n}^*\},
\]
where
$p_{S,S'}$
is the projection map related to the direct sum $S\oplus S'$.
\end{thm}

The polynomial representation of the projection is of type
\[
p_{S,S'}=\sum_{j=0}^{n-1}\left(\sum_{h=s}^{n-1}\xi_h{\xi_h^*}^{q^j}\right)X^{q^j},
\]
where $\{\xi_1,\xi_2,\ldots,\xi_n\}$ and $\{\xi_1^*,\xi_2^*,\ldots,\xi_n^*\}$ are suitable dual $\Fq$-bases of $\Fqn$ related to $S$ and $S'$.
We are seeking shorter polynomial representations for specific linear sets with complementary weights.

\begin{thm}\cite[Theorem 4.5]{NapPolSanZul22}
Let  $1 < t < n$ and  $n = \ell t$. There exist $\Fq$-linear sets of rank $\rho$ in $\PG ( 1 , q^n )$ with one
point of weight $t$, one point of weight $s$ and all others of weight one for the following values of $n$, $\rho$ and $s$:
\begin{itemize}
\item $n$ even, $\rho = t + s$ and any $s\in\{ 1 ,\ldots , n/2 \}$;
\item $n$ odd, $\rho = t + s$ and any $s \in \{ 1 ,\ldots , \frac{n-t}2\}$.
\end{itemize}

\end{thm}

Previous result provides the following consequence on regular fat linear sets.

\begin{cor}\label{c:2t}
If $t$ divides $n$ and $t\neq n$, then there is a $(2,t)$-regular fat linear set of rank $2t$ in $\PG(1,q^n)$.
\end{cor}

In Subsection~\ref{ss:construction1} we will find an explicit polynomial representation for some of the $(2,t)$-regular fat linear sets mentioned in Corollary~\ref{c:2t}.

In \cite[Theorem 3.5]{Zul23} examples of $(k,i)$-regular fat linear sets of rank $ki$ with $i\le t$ in $\PG(k-1,q^{2t})$ are described.


\subsection{Further constructions of \texorpdfstring{$(2,i)$}{(2,i)}-regular fat linear sets}

In the context of representability of the direct sum of uniform $q$-matroids, \cite[Theorem 4.4]{alfarano} provides examples of $\fq$-linear sets in $\mathrm{PG}(1,q^n)$ with only two points of weight greater than one, which have weights $i_1$ and $i_2$, of rank $i_1+i_2$. We omit details on the relationship between $q$-matroids and regular fat linear sets for the sake of brevity. We state \cite[Theorem 4.4]{alfarano} in the case $i_1 = i_2$, i.e., when \cite[Theorem 4.4]{alfarano} provides examples of $(2,i)$-regular fat linear sets.


\begin{thm}
There exist $(2,i)$-regular fat linear sets in $\PG(1,q^n)$ in the following cases:
\begin{enumerate}
       \item $n$  even and $n\geq 2i$; 
       \item $n\geq i^2$;
       \item $n=t_1t_2$ with $t_1\ge i$ and $i\le 1+t_1(t_2-1)/2$;
       \item $n=t_1t_2$ with $t_1\geq i$ and $t_2\geq 2$;
       \item $q=p^h$, $n=p^r$, $2i-1\le n/2$.
\end{enumerate}
\end{thm}

To our knowledge, there are no explicit examples of $(r,i)$-regular fat linear sets in $\PG(1,q^n)$ with $r>2$ and $i>2$ in the literature.
In the next section, we will present constructions of these linear sets.

\section{\texorpdfstring{Examples of $(r,i)$-regular fat linear sets for $r,i>2$}{Examples of (r,i)-regular fat linear sets for r,i>2}}\label{s:examples}
In this section, we present two new examples of regular fat linear sets. Subsections \ref{ss:construction1} and \ref{ss:construction2}, and more specifically Theorems \ref{t:mainexample1} and \ref{t:mainexample2}, are devoted to constructions in $\PG(k-1,q^{\ell t})$ for $\ell=2$ and $\ell>2$, respectively.
For future reference, we quote the following results by Gow and Quinlan.
\begin{thm}\cite{GQ1,GQ2}\label{t:GQ}
Let $q$ be a power of a prime, $s$ and $n$ be two coprime positive integers, and $f=\sum_{j=0}^dc_jX^{q^{js}}\in\Fqn[X]$, $c_d\neq0$. Then
\begin{enumerate}[$(i)$]
\item $f$ has at most $q^d$ roots in $\Fqn$;
\item if $f$ has $q^d$ roots in $\Fqn$, then $\N_{q^n/q}(c_d)=(-1)^{nd}\N_{q^n/q}(c_0)$.
\end{enumerate}
\end{thm}

\subsection{A construction in \texorpdfstring{$\PG(k-1,q^{2t})$}{PG(k-1,q2t)}}\label{ss:construction1}

\begin{thm}\label{t:mainexample1}
Let $q$ be odd, $t\ge3$, $\gcd(s,t)=1$,\\ $w\in E=\{x\in\F_{q^{2t}}\colon\Tr_{q^{2t}/q^t}(x)=0\}$, $w\neq0$, $\operatorname N_{q^t/q}(w^2)\neq(-1)^t$, and let $I$ be an $i$-dimensional $\Fq$-subspace of $\F_{q^t}$, $i\geq2$. Define
\begin{equation}\label{e:defT} T=T_{s,w,I}=\{x+wx^{q^s}\colon x\in I\}, \end{equation}
then for any $k>1$, $L_{T^k}$ is a $((q^k-1)/(q-1),i)$-regular fat linear set of rank $ki$ in $\PG(k-1,q^{2t})$.
The points of weight $i$ are precisely the elements of $\PG(k-1,q)$, i.e., $(a_1:a_2:\ldots:a_k)$ with $(0,0,\ldots,0)\neq(a_1,a_2,\ldots,a_k)\in\F_q^k$.
\end{thm}
\begin{proof}
The assertion for $\dim_{\Fq}(T^k)=ki$ is clear.

\noindent
1) We prove that if $P=(a_1:a_2:\ldots:a_k)$ and $(0,0,\ldots,0)\neq(a_1,a_2,\ldots,a_k)\in\F_q^k$, then $w_{L_{T^k}}(P)=i$.
The intersection $P\cap T^k$ can be derived from the simultaneous equations
\begin{equation}\label{e:sistperw}
(\rho'+w\rho'')a_j=x_j+wx^{q^s}_j,\quad j=1,\ldots,k,
\end{equation}
in the unknowns $\rho',\rho''\in\Fqt$, $x_j\in I$.
For all $j$ such that $a_j\neq0$ such equations are equivalent to
\[
\begin{cases}\rho'a_j&=x_j\\ \rho''a_j&=x_j^{q^s}\end{cases}
\quad\Leftrightarrow\quad
\begin{cases}\rho''&=\rho'^{q^s}\\ x_j&=\rho' a_j.\end{cases}
\]
Therefore,
\[
P\cap T^k=\left\{(\rho'+w\rho'^{ q^s})(a_1,a_2,\ldots,a_k)\colon\rho'\in I\right\}.
\]
Hence $w_{L_{T^k}}(P)=i$.

\noindent
2) Assume now that $Q=(b_1:b_2:\ldots:b_k)\notin\PG(k-1,q)$.
It may be assumed that $b_h=1$, $b_r\notin\Fq$ for some $h,r\in\{1,\ldots,k\}$. Let $b_r=\alpha'+w\alpha''$, $\alpha',\alpha''\in\Fqt$.
From the simultaneous equation similar to \eqref{e:sistperw}, taking only $j=h$ and $j=r$,
\[
\begin{cases}\rho'&=x_h\\ \rho''&=x_h^{q^s}\\
\rho'\alpha'+w^2\rho''\alpha''&=x_r\\
\rho'\alpha''+\rho''\alpha'&=x_r^{q^s}\end{cases}
\quad\Leftrightarrow\quad
\begin{cases}\rho'&=x_h\\ \rho''&=x_h^{q^s}\\
\alpha'x_h+w^2\alpha''x_h^{q^s}&=x_r\\
\alpha''x_h+\alpha'x_h^{q^s}&=x_r^{q^s}.\end{cases}
\]
Combining the last two equations, we obtain the following equation in $x_h$, whose coefficients are in $\Fqt$:
\[
w^{2q^s}\alpha''^{q^s}x_h^{q^{2s}}+(\alpha'^{q^s}-\alpha')x_h^{q^s}-\alpha''x_h=0.
\]
If $w_{L_{T^k}}(Q)>1$, then this equation has more than $q$ solutions in $\Fqt$. The condition $b_r\notin\Fq$ implies that not all coefficients are zero.
By Theorem~\ref{t:GQ}, a necessary condition for having more than $q$ solutions, that is, $q^2$ solutions, is
$\N_{q^t/q}(w^{2q^s}\alpha''^{q^s})=\N_{q^t/q}(-\alpha'')$, $\alpha''\neq0$, equivalent to $\N_{q^t/q}(w^2)=(-1)^t$ and this is against the hypotheses.
\end{proof}

Our next goal is to prove that the assumptions in the previous theorem can be fulfilled.

\begin{prop}\label{p:haresatisfied}
    Let $q$ be odd.
    For $q>3$ there exists $w\in E$, $w\neq0$ such that $\operatorname N_{q^t/q}(w^2)\neq-1$, while $\operatorname N_{q^t/q}(w^2)\neq1$ holds for any $w\in E$.
\end{prop}
\begin{proof}
    Let $\omega$ be a generator of the multiplicative group of $\F_{q^{2t}}$, and define $w=\omega^{(q^t+1)/2}$.
    We have $w^{q^t-1}=\omega^{(q^{2t}-1)/2}=-1$, and hence $w\in E$.
    Furthermore,
    \[
    \N_{q^t/q}(w^2)=\omega^{\frac{q^{2t}-1}{q-1}},
    \]
    and this is equal to $-1$ if and only if $q=3$.
Now, let $w$ be any nonzero element of $E$.
Since the elements in $\F_{q^t}^*$ are of type $\omega^{j(q^t+1)}$, $j\in\mathbb Z$, we have
\[
w=\omega^{\frac{(2j+1)(q^t+1)}2}.
\]
Then $\N_{q^t/q}(w^2)=\omega^{(q^{2t}-1)(2j+1)/(q-1)}$ and this equals one if and only if $q-1$, which is even, divides $2j+1$. So, $\N_{q^t/q}(w^2)\neq1$.
\end{proof}

\begin{prop}\label{prop:rqt}
    The linear sets $L_{T^k}$ in Theorem~\ref{t:mainexample1} are $\Rt$partially scattered.
\end{prop}
\begin{proof}
First, we require to find the size of the intersection
\[  T^k\cap  \la(a_1,a_2,\ldots,a_k)\ra_{\Fqt},\]
for all $a_1$, $a_2$, $\ldots$, $a_k\in\F_{q^{2t}}$.
This gives $k$ equations $\theta a_j=x_j+wx_j^{q^s}$, $\theta\in\Fqt$, $x_j\in I$, $j=1,\ldots,k$.
Take $a_j\neq0$. Let $a_j=a'+wa''$, $a',a''\in\Fqt$.
It results in $\theta a'=x_j$, $\theta a''=x_j^{q^s}$, hence $\theta^{q^s}a'^{q^s}-\theta a''=0$ and this has  at most $q$ solutions in $\theta$.

\end{proof}

\begin{thm}\label{t:polform}
For any $\mu\in\F_{q^t}$ such that $\N_{q^t/q}(\mu)=1$, $\mu\neq1$, the linear set $L_{T^2}\subseteq\PG(1,q^{2t})$ with $T=T_{s,w,\Fqt}$ is equivalent up to the action of $\operatorname{GL}(2,q^{2t})$ to $L_f$, where
\begin{align}
f=&(\mu^{q^s}-1)\left((\mu+1)X^{q^t}-2w^{-q^{t-s}}(X^{q^{t-s}}-X^{q^{2t-s}})\right)\notag\\
&+(\mu-1)\left((\mu^{q^s}+1)X^{q^t}+2w\mu^{q^s}(X^{q^s}+X^{q^{t+s}})\right).\label{e:pol1}
\end{align}
For even $t$, an admissible choice for $f$ is 
\begin{equation}\label{e:pol2}
f=X^{q^s}+X^{q^{t+s}}+m(X^{q^{t-s}}-X^{q^{2t-s}}),
\end{equation}
where $m$ is the $(q+1)$-th power of a nonzero element of $E$.
\end{thm}
\begin{proof}
    Let $\eta\in\Fqt$ such that $\eta^{q^s-1}=\mu^{q^s}$.
    By substitution $v=\eta z$,
    \begin{eqnarray*}
    T^2&=&\left\{(u+wu^{q^s},v+wv^{q^s})\colon u,v\in\Fqt\right\}\\
    &=&\left\{(u+wu^{q^s},\eta(z+w\mu^{q^s}z^{q^s}))\colon u,z\in\Fqt\right\},
    \end{eqnarray*}
    that is GL-equivalent to
    $\left\{(u+wu^{q^s},z+w\mu^{q^s}z^{q^s})\colon u,z\in\Fqt\right\}$.
After a left multiplication by $\begin{pmatrix}1&1\\ 1&-1\end{pmatrix}$ we obtain the following:
\[
\left\{\left(u+z+w(u+\mu z)^{q^s},u-z+w(u-\mu z)^{q^s}\right)\colon u,z\in\Fqt\right\}.
\]
Now, let $x',x''\in\Fqt$.
Equation $u+z+w(u+\mu z)^{q^s}=x'+wx''$, $u,z\in\Fqt$, has the unique solution
\[
\begin{cases}
u&=(\mu x'-{x''}^{q^{t-s}})/(\mu-1)\\
z&=({x''}^{q^{t-s}}-x')/(\mu-1).
\end{cases}
\]
On the other hand, equation $x'+wx''=x$, $x\in\F_{q^{2t}}$, where $x'$ and $x''$ are unknown elements of $\Fqt$, has as unique solution
\[
\begin{cases}
x'&=(x+x^{q^t})/2\\
x''&=w^{-1}(x-x^{q^t})/2.
\end{cases}
\]
Therefore, $x=u+z+w(u+\mu z)^{q^s}$ is equivalent to
\begin{eqnarray*}u&=\frac1{2(\mu-1)}\left(\mu(x+x^{q^t})-w^{-q^{t-s}}(x^{q^{t-s}}-x^{q^{2t-s}})\right),\\
z&=\frac1{2(\mu-1)}\left(w^{-q^{t-s}}(x^{q^{t-s}}-x^{q^{2t-s}})-(x+x^{q^{t}})\right).
\end{eqnarray*}
This leads to the representation of the linear set as $L_g$ where
\begin{eqnarray*}
g(x)&=&u-z+w(u^{q^s}-\mu^{q^s}z^{q^s})\\
&=&\frac1{2(\mu-1)}\left((\mu+1)(x+x^{q^t})-2w^{-q^{t-s}}(x^{q^{t-s}}-x^{q^{2t-s}})\right)\\
&&+\frac w{2(\mu^{q^s}-1)}\left(2\mu^{q^s}(x^{q^s}+x^{q^{t+s}})-(\mu^{q^s}+1)w^{-q^t}(x^{q^t}-x)\right).
\end{eqnarray*}
The monomial in $x$ may be neglected.
Taking into account $w^{-q^t}=-w^{-1}$, this leads to \eqref{e:pol1}.
For $t$ even and $\mu=-1$, multiplying by $2w^{-1}$ and setting $m=w^{-(q^{t-s}+1)}$ leads to \eqref{e:pol2}.
\end{proof}

Now, we will address the question of equivalence between the linear sets that we constructed in Theorem~\ref{t:mainexample1}.

\begin{thm}
    Assume that $s$ and $\tilde s$ are coprime with $t$, $1\le s<t$, $1\le\tilde s<t$, $w,\tilde w\in E\setminus\{0\}$.
    Let $I$, $\tilde I$ be $\Fq$-subspaces of $\Fqt$ with $\dim_{\Fq}I>2$. Then
    \begin{enumerate}[$(i)$]
    \item If $s\neq\tilde s,t-\tilde s$, and $I=\Fqt$, then $\left(T_{s,w,I}\right)^k$ and $\left(T_{\tilde s,\tilde w,\tilde I}\right)^k$ belong to distinct orbits under the action of $\operatorname{\Gamma L}(k,q^{2t})$.
        \item The $\Fq$-subspaces $\left(T_{s,w,I}\right)^k$ and $\left(T_{s,\tilde w,\tilde I}\right)^k$ are equivalent under $\operatorname{\Gamma L}(k,q^{2t})$ if, and only if, there is an automorphism $\iota$ of $\F_{q^{2t}}$ such that $\N_{q^t/q}(w^\iota\tilde w^{-1})=1$, $\tilde I=\Theta I^\iota$ where $\Theta\in\Fqt$ and $\Theta^{q^s-1}=w^\iota\tilde w^{-1}$.
    \item The $\Fq$-subspaces  $\left(T_{s,w,I}\right)^k$ and $\left(T_{t-s,\tilde w,\tilde I}\right)^k$ are equivalent under $\operatorname{\Gamma L}(k,q^{2t})$ if, and only if, there is an automorphism $\iota$ of $\F_{q^{2t}}$ such that $\N_{q^t/q}(w\ii\tilde w)=(-1)^{t-s}$, $\tilde I=w\ii w_0 {I\ii}^{q^s}$ where $w_0\in E$ and $w_0^{1-q^{t-s}}={w^\iota}^{q^{t- s}}\tilde w$.
    \end{enumerate}
\end{thm}
\begin{proof}
    Assume that $\left(T_{s,w,I}\right)^k$ and $\left(T_{\tilde s,\tilde w,\tilde I}\right)^k$ are equivalent under $\operatorname{\Gamma L}(k,q^{2t})$, that is, there are a matrix $A\in\GL(k,q^{2t})$ and an automorphism $\iota$ of $\F_{q^{2t}}$ such that for any $u_1,\ldots,u_k\in I$ there exist $x_1,\ldots,x_k\in\tilde I$ satisfying
    \begin{equation}\label{e:sonoequiv}
    A
    \begin{pmatrix}u_1\ii+w\ii {u_1\ii}^{q^s}\\ u_2\ii+w\ii {u_2\ii}^{q^s}\\ \vdots\\ u_k\ii+w\ii {u_k\ii}^{q^s}\end{pmatrix}=\begin{pmatrix}x_1+\tilde w x_1^{q^{\tilde s}}\\ x_2+\tilde w x_2^{q^{\tilde s}}\\ \vdots\\ x_k+\tilde w x_k^{q^{\tilde s}}\end{pmatrix},\end{equation} and conversely.
Let $a_{jh}$ denote the element of $A$ in the $j$-th row and the $h$-th column.
Let $j,h\in\{1,\ldots,k\}$, and $a_0\in\Fqt$, $a_1\in E$ be such that $a_{jh}=a_0+a_1$.
Since \eqref{e:sonoequiv} also holds for $u_1=u_2=\ldots=u_{h-1}=u_{h+1}=\ldots=u_k=0$, for any $u=u_h\in I$ there exists $x=x_j\in\tilde I$ satisfying
\begin{equation}\label{e:v01}
\begin{cases}x&=a_0u\ii+w\ii a_1{u\ii}^{q^s}\\ \tilde wx^{q^{\tilde s}}&=a_1u\ii+w\ii a_0{u\ii}^{q^s}.\end{cases}
\end{equation}
From \eqref{e:v01},
\begin{equation}\label{e:v03}
\tilde w{w\ii}^{q^{\tilde s}}a_1^{q^{\tilde s}}{u\ii}^{q^{s+\tilde s}}+\tilde wa_0^{q^{\tilde s}}{u\ii}^{q^{\tilde s}}-w\ii a_0{u\ii}^{q^s}-a_1u\ii=0,
\end{equation}
for all $u\in I$.
If $s\neq\tilde s,t-\tilde s$, taking into account the assumption $I=\Fqt$ this implies $a_{jh}=0$ for all $j$ and $h$. This proves~$(i)$.

Now we prove $(ii)$.
Assume that $\left(T_{s,w,I}\right)^k$ and $\left(T_{s,\tilde w,\tilde I}\right)^k$ are equivalent under $\operatorname{\Gamma L}(k,q^{2t})$.
Starting from \eqref{e:sonoequiv} with $s=\tilde s$, by \eqref{e:v01}, for all $u\in I$,
\[
\tilde w {w\ii}^{q^s}a_1^{q^s}{u\ii}^{q^{2s}}+(\tilde wa_0^{q^s}-w\ii a_0){u\ii}^{q^s}-a_1 u\ii =0.
\]
By Theorem~\ref{t:GQ} and $\dim_{\Fq}I>2$, it follows that for any $j,h\in\{1,\ldots,k\}$, $a_{jh}\in\Fqt$ and $a_{jh}(a_{jh}^{q^s-1}-w\ii\tilde w^{-1})=0$.
Since $A$ has some non-zero entries, this implies that $\N_{q^t/q}(w\ii\tilde w^{-1})=1$.
Furthermore, $a_{jh}\in\la \Theta\ra_{\Fq}$,  where $\Theta^{q^s-1}=w^\iota\tilde w^{-1}$, for all $j$ and $h$.
By \eqref{e:v01},
$ a_{jh}I\ii\subseteq \tilde I$
and the equality must hold for the non-zero values among the $a_{jh}$s.

Conversely, assume that the conditions in $(ii)$ hold.
The equation
\[
\Theta \begin{pmatrix}u_1\ii+w\ii {u_1\ii}^{q^s}\\ u_2\ii+w\ii {u_2\ii}^{q^s}\\ \vdots\\u_k\ii+w\ii {u_k\ii}^{q^s}\end{pmatrix}=\begin{pmatrix}x_1+\tilde w x_1^{q^s}\\ x_2+\tilde w x_2^{q^s}\\ \vdots\\ x_k+\tilde w x_k^{q^s}\end{pmatrix}
,\]
is equivalent to
\[
\begin{cases}
    x_j&=\Theta u_j\ii\\ \tilde w \Theta^{q^s}{u_j\ii}^{q^s}&=w\ii \Theta{u_j\ii}^{q^s},
\end{cases}\quad j=1,\ldots,k.
\]
Since $\tilde w \Theta^{q^s}=w\ii \Theta$ and $\Theta u_j\ii\in\tilde I$, $j=1,\ldots,k$, the semilinear map $(\delta_1,\delta_2,\ldots,\delta_k)\mapsto \Theta(\delta_1\ii,\delta_2\ii,\ldots,\delta_k\ii)$ transforms $\left(T_{s,w,I}\right)^k$ into $\left(T_{s,\tilde w,\tilde I}\right)^k$.

The proof of assertion $(iii)$ is similar. 
From
\[
    A
    \begin{pmatrix}u_1\ii+w\ii {u_1\ii}^{q^s}\\ u_2\ii+w\ii {u_2\ii}^{q^s}\\ \vdots\\ u_k\ii+w\ii {u_k\ii}^{q^s}\end{pmatrix}=\begin{pmatrix}x_1+\tilde w x_1^{q^{t- s}}\\ x_2+\tilde w x_2^{q^{t- s}}\\ \vdots\\ x_k+\tilde w x_k^{q^{t- s}}\end{pmatrix},
\]
we obtain
\begin{equation}\label{e:v04}
\begin{cases}x&=a_0u\ii+w\ii a_1{u\ii}^{q^s}\\ \tilde wx^{q^{t- s}}&=a_1u\ii+w\ii a_0{u\ii}^{q^s}.\end{cases} 
\end{equation}
By \eqref{e:v04}, all $u\in I$ satisfy the following equation:
\begin{equation}\label{e:v05}
(\tilde w{w\ii}^{q^{t-s}}a_1^{q^{t-s}}-a_1)u\ii+\tilde wa_0^{q^{t-s}}{u\ii}^{q^{t-s}}-w\ii a_0{u\ii}^{q^s}=0.
\end{equation}
Then the following equation, obtained from \eqref{e:v05} by making $z={u\ii}^{q^{t-s}}$, has more than $q^2$ solutions:
\[
-w\ii a_0z^{q^{2s}}+(\tilde w{w\ii}^{q^{t-s}}a_1^{q^{t-s}}-a_1)z^{q^s}+\tilde wa_0^{q^{t-s}}z=0.
\]
Therefore, for all $j,h=1,\ldots,k$, we have $a_{jh}\in E$ and $a_{jh}(\tilde w{w\ii}^{q^{t-s}}a_{jh}^{q^{t-s}-1}-1)=0$. This implies $a_{jh}\in\la w_0\ra_{\Fq}$ where $w_0\in E$ and $w_0^{1-q^{t-s}}={w^\iota}^{q^{t- s}}\tilde w$. The first equation in \eqref{e:v04} yields $w\ii w_0{I\ii}^{q^s}=\tilde I$.
 Setting $\tau=q^{t-s}$,
\[
\N_{q^t/q}(w\ii \tilde w)=\N_{q^t/q}({w\ii}^\tau \tilde w)=\N_{q^t/q}(a_{jh}^{1-\tau})=a_{jh}^{(1-\tau)\frac{\tau^t-1}{\tau-1}}=a_{jh}^{1-\tau^t}.
\]
Taking into account
\[ a_{ij}^{\tau^t}=a_{ij}^{({q^t})^{t-s}}=(\ldots((a_{jh}^{q^t})^{q^t})\ldots)^{q^t}=(-1)^{t-s}a_{jh},\]
we obtain $\N_{q^t/q}(w\ii \tilde w)=(-1)^{t-s}$.

Conversely, assume that the conditions in $(iii)$ hold.
The equation
\[
w_0\begin{pmatrix}u_1\ii+w\ii {u_1\ii}^{q^s}\\ u_2\ii+w\ii {u_2\ii}^{q^s}\\ \vdots\\ u_k\ii+w\ii {u_k\ii}^{q^s}\end{pmatrix}=\begin{pmatrix}x_1+\tilde w x_1^{q^{t-s}}\\ x_2+\tilde w x_2^{q^{t-s}}\\ \vdots\\ x_k+\tilde w x_k^{q^{t-s}}\end{pmatrix},
\]
is equivalent to
\[
\begin{cases}
    x_j&=w_0w\ii{u_j\ii}^{q^s}\\ \tilde w x_j^{q^{t-s}}&=w_0 u_j\ii,
\end{cases}\quad j=1,\ldots,k.
\]
The second equation is equivalent to the first.
Therefore, the semilinear map $(\alpha_1,\alpha_2,\ldots,\alpha_k)\mapsto w_0(\alpha_1\ii,\alpha_2\ii,\ldots,\alpha_k\ii)$ transforms $\left(T_{s,w,I}\right)^k$ into $\left(T_{t-s,\tilde w,\tilde I}\right)^k$.
\end{proof}
\begin{rem}
    Regarding the assertion $(iii)$ in the previous proposition, note that if $y\in\Fqt$ and $\N_{q^t/q}(y)=(-1)^{t-s}$, then a $w_0\in E$ exists such that $w_0^{1-q^{t-s}}=y$.
    In fact, the norm $\N_{q^t/q}(w^{1-q^{t-s}})$ is equal to $(-1)^{t-s}$ for any $w\in E$, $w\neq0$, and both sets $A=\{w^{1-q^{t-s}}\colon w\in E,\,w\neq0\}$ and $B=\{y\in\Fqt\colon\N_{q^t/q}(y)=(-1)^{t-s}\}$ have size $(q^t-1)/(q-1)$. Hence, $A=B$.
\end{rem}

\subsection{A construction in \texorpdfstring{$\PG(k-1,q^{\ell t})$}{PG(k-1,q^{lt})}, \texorpdfstring{$\ell>2$}{l>2}}\label{ss:construction2}

Let $\ell$ and $t$ be two positive integers and assume that $\ell$ divides $q^t-1$.
Fix $\eps\in\Fqt$ such that $\{1,\eps,\eps^2,\ldots,\eps^{\ell-1}\}$ is the set of all roots of $X^\ell-1$ and define
\begin{equation}\label{e:ej}E_j=\{x\in\F_{q^{\ell t}}\colon x^{q^t}-\eps^jx=0\},\quad j=0,\ldots,\ell-1.\end{equation}
Since $\N_{{q^{\ell t}}/q^t}(\eps^j)=1$ for all $j$, each $E_j$ is a one-dimensional $\Fqt$-subspace of $\F_{q^{\ell t}}$. Hence, there exists $\eta\in E_1$ such that 
\begin{equation}\label{e:E_i}
E_j=\langle\eta^j\rangle_{\Fqt}, \quad \,\,j=0,\ldots,\ell-1.
\end{equation}

\begin{lem}\begin{enumerate}[$(i)$]
\item For all $j$ and $h$, $E_{j+h}=E_jE_h$, where the indices are taken $\Mod\ell$.
\item \[\mathbb{F}_{q^{\ell t}}=E_0\oplus E_1\oplus\ldots\oplus E_{\ell-1}.\]
\end{enumerate}
\end{lem}
\begin{proof} Since assertion $(i)$ is clear, we only prove $(ii)$.
    Taking into account the dimensions of $\mathbb F_{q^{\ell t}},E_0, \ldots,E_{\ell-1}$, it suffices to prove, by induction on $j=0,1,\ldots,\ell-1$, that the sum $E_0+E_1+\cdots+E_j$ is direct. Assume that $E_0+ E_1+\ldots+ E_j$ is a direct sum for $j\in\{0,\ldots,\ell-2\}$. We claim $E_{j+1}\cap (E_0+ E_1+\ldots + E_j)=\{0\}$. Consider $E_{j+1}\ni x=x_0+x_1+\ldots+x_j$, where $x_h\in E_h$, $h=0,\ldots,j$. By raising to the $q^t$-th power, we get $x^{q^t}=x_0^{q^t}+x_1^{q^t}+\ldots+x_j^{q^t}$; that is,
    \[\varepsilon^{j+1}x=x_0+\varepsilon x_1+\ldots+\varepsilon^j x_j.\]
    Subtracting from the latter equation the identity $\varepsilon^{j+1}x=\varepsilon^{j+1}(x_0+x_1+\ldots+x_j)$, we get
    \[0=(1-\varepsilon^{j+1})x_0+(\varepsilon-\varepsilon^{j+1})x_1+\ldots+(\varepsilon^j-\varepsilon^{j+1})x_j.\]
    Since $x_h(\varepsilon^h-\varepsilon^{j+1})\in E_h$ for all $h=0,\ldots,j$, and since by the induction hypothesis $E_0+ E_1+\ldots+ E_{j}$ is a direct sum, we get the thesis as $x_0=x_1=\ldots=x_j=0$.
\end{proof}

We introduce a new class of regular fat linear sets.

\begin{thm}\label{t:mainexample2}
Let  $\gcd(s,t)=1$, $\ell>2$, $\ell \mid q^t-1$, $\eta$ be defined as in Equation \eqref{e:E_i}, and $I$ be an $i$-dimensional $\Fq$-subspace of $\F_{q^t}$, $i\geq2$. Define
\begin{equation}\label{e:defT2} T=T_{s,\eta,I}=\{x+\eta x^{q^s}\colon x\in I\}, \end{equation}
then for any $k>1$, $L_{T^k}$ is a $((q^k-1)/(q-1),i)$-regular fat linear set of rank $ki$ in $\PG(k-1,q^{\ell t})$.
The points of weight $i$ are precisely the elements of $\PG(k-1,q)$, i.e.\ $(a_1:a_2:\ldots:a_k)$ with $(0,\ldots,0)\neq(a_1,\ldots,a_k)\in\F_q^k$.
\end{thm}

\begin{proof}
    Similar as in the previous construction the assertion $\dim_{\Fq}(T^k)=ki$ is clear.

\noindent
We prove that if $P=(a_1:a_2:\ldots:a_k)$ and $(0,\ldots,0)\neq(a_1,\ldots,a_k)\in\F_q^k$, then $w_{L_{T^k}}(P)=i$.
The intersection $P\cap T^k$ can be derived from the simultaneous equations
\begin{equation}\label{e:sistperw2}
\sum_{h=0}^{\ell-1}\beta_h a_j=x_j+\eta x^{q^s}_j,\quad j=1,\ldots,k,
\end{equation}
in the unknowns $\beta_h\in E_h$, $x_j\in I$.
For all $j$ such that $a_j\neq0$ such equations are equivalent to
\[
\begin{cases}\beta_0a_j=x_j\\ \beta_1a_j=\eta x_j^{q^s}\\ \beta_2=\ldots=\beta_{\ell-1}=0\end{cases}
\quad\Leftrightarrow\quad
\begin{cases}\beta_1=\eta\beta_0^{q^s}\\ x_j=\beta_0 a_j\\ \beta_2=\ldots=\beta_{\ell-1}=0.\end{cases}
\]
Therefore,
\[
P\cap T^k=\left\{(\beta_0+\eta\beta_0^{ q^s})(a_1,a_2,\ldots,a_k)\colon\beta_0\in I\right\}.
\]
Hence, $w_{L_{T^k}}(P)=i$.

\noindent
Assume now that $(b_1:b_2:\ldots:b_k)\notin\PG(k-1,q)$.
It may be assumed w.l.o.g.\ that $b_h=1$, $b_r\notin\Fq$ for some $h,r\in\{1,\ldots,k\}$. Let $b_r=\sum_{j=0}^{\ell-1}z_j$, $z_j\in E_j$, $j=1,\ldots,\ell-1$.
From the simultaneous equations similar to \eqref{e:sistperw2}, taking only $j=h$ and $j=r$, we have the following equations in the unknowns $\beta\in\F_{q^{\ell t}}$, $x_h,x_r\in I$:
\begin{equation}\label{projections}
\begin{cases}\beta b_h&=x_h+\eta x_h^{q^s}\\ \beta b_r&=x_r+\eta x_r^{q^s}\end{cases}
\quad\Leftrightarrow\quad
\begin{cases}x_h+\eta x_h^{q^s}=&\beta\\
x_r+\eta x_r^{q^s}=&\sum_{j=0}^{\ell-1}z_j(x_h+\eta x_h^{q^s}).\end{cases}
\end{equation}
The projection on $E_2$ of the second equation of \eqref{projections} reads
\[0=z_2x_h+z_1\eta x_h^{q^s}.\]
If $z_1\neq0$ or $z_2\neq0$ the equation has at most $q$ solutions, giving the thesis. On the other hand, if $z_1=z_2=0$, the projection on $E_1$ of the second equation in~\eqref{projections} reads
\[\eta x_r^{q^s}=z_0 \eta x_h^{q^s},\]
implying $x_r=z_0^{q^{t-s}}x_h$. By replacing this value in the second equation in~\eqref{projections} we get the projection on $E_0=\Fqt$ 
\[z_0^{q^{t-s}}x_h=z_0 x_h+z_{\ell-1}\eta x_h^{q^s}.\]
If we suppose to have more than $q$ solutions, we find $z_0\in\fq$ and $z_{\ell-1}=0$, hence $x_r\in\langle x_h\rangle_{\fq}$, contradicting $b_r\notin\fq$.
\end{proof}

\begin{prop}
    The linear sets $L_{T^k}$ in Theorem~\ref{t:mainexample2} are $\Rt$partially scattered.
\end{prop}

\begin{proof}
    First, we require to find the size of \[T^k\cap\la(a_1,a_2,\ldots,a_k)\ra_{\Fqt}\] for all $a_1$, $a_2$, $\ldots$, $a_k\in\F_{q^{2t}}$.
    This gives equations
    \begin{equation}\label{e:ottobre}x_j+\eta x_j^{q^s}=\theta a_j\end{equation} in unknowns $x_j\in I$, $j=1,\ldots,k$ and $\theta\in\Fqt$.
    Take $a_j\neq0$, $a_j=\sum_{h=0}^{\ell-1}z_h$, with $z_h\in E_h$ for $h=0,\ldots,\ell-1$.
    If $z_0=z_1=0$, then \eqref{e:ottobre} has a unique solution $\theta=0$. Otherwise, from $x_j=\theta z_0$, $\eta x_j^{q^s}=\theta z_1$ we have $\eta \theta^{q^s}z_0^{q^s}-\theta z_1=0$ which is a nonvanishing equation in $\theta$ and has at most $q$ solutions.
\end{proof}

Now, we will address the question of equivalence between the linear sets that we constructed in Theorem~\ref{t:mainexample2}.

\begin{thm}
    Assume that $s$ and $\tilde s$ are coprime with $t$, $1\le s<t$, $1\le\tilde s<t$, $\eta,\tilde\eta\in E_1\setminus\{0\}$.
    Let $I$, $\tilde I$ be $\Fq$-subspaces of $\Fqt$ with $\dim_{\Fq}I>2$. Then the following properties hold:
    \begin{enumerate}[$(i)$]
    \item If $s\neq\tilde s,t-\tilde s$, and $I=\Fqt$, then $\left(T_{s,\eta,I}\right)^k$ and $\left(T_{\tilde s,\tilde\eta,\tilde I}\right)^k$ belong to distinct orbits under the action of $\operatorname{\Gamma L}(k,q^{\ell t})$.
        \item The $\Fq$-subspaces $\left(T_{s,\eta,I}\right)^k$ and $\left(T_{s,\tilde\eta,\tilde I}\right)^k$ are equivalent under $\operatorname{\Gamma L}(k,q^{\ell t})$ if, and only if, there is an automorphism $\iota$ of $\F_{q^{\ell t}}$ such that $E_1^{\iota}=E_1$, $\N_{q^t/q}(\eta^\iota\tilde\eta^{-1})=1$, $\tilde I=\Theta I^\iota$ where $\Theta\in\Fqt$ and $\Theta^{q^s-1}=\eta^\iota\tilde\eta^{-1}$.
    \item The $\Fq$-subspaces  $\left(T_{s,\eta,I}\right)^k$ and $\left(T_{t-s,\tilde\eta,\tilde I}\right)^k$ are equivalent under $\operatorname{\Gamma L}(k,q^{\ell t})$ if, and only if, there is an automorphism $\iota$ of $\F_{q^{\ell t}}$ such that $E_1^{\iota}=E_{\ell-1}$, $\tilde I=\Lambda\eta\ii {I\ii}^{q^s}$ where $\Lambda\in E_1$ and $\Lambda^{1-q^{t-s}}=\tilde\eta{\eta^\iota}^{q^{t- s}}$.
    \end{enumerate}
\end{thm}
\begin{proof}
    Assume that $\left(T_{s,\eta,I}\right)^k$ and $\left(T_{\tilde s,\tilde\eta,\tilde I}\right)^k$ are equivalent under $\operatorname{\Gamma L}(k,q^{\ell t})$, that is, there are a matrix $A\in\GL(k,q^{\ell t})$ and an automorphism $\iota$ of $\F_{q^{\ell t}}$ such that for any $u_1,\ldots,u_k\in I$ there exist $x_1,\ldots,x_k\in\tilde I$ satisfying
    \begin{equation}\label{e:sonoequiv2}
    A
    \begin{pmatrix}u_1\ii+\eta\ii {u_1\ii}^{q^s}\\ u_2\ii+\eta\ii {u_2\ii}^{q^s}\\ \vdots\\ u_k\ii+\eta\ii {u_k\ii}^{q^s}\end{pmatrix}=\begin{pmatrix}x_1+\tilde\eta x_1^{q^{\tilde s}}\\ x_2+\tilde\eta x_2^{q^{\tilde s}}\\ \vdots\\ x_k+\tilde\eta x_k^{q^{\tilde s}}\end{pmatrix},\end{equation}
    and conversely. Let $a_{jh}$ denote the element of $A$ in the $j$-th row and the $h$-th column.
Let $\alpha_\nu\in E_\nu$, $\nu=0,\ldots,\ell-1$ such that $a_{jh}=\sum_{\nu=0}^{\ell-1}\alpha_\nu$. Note that the automorphism $\iota$ fixes $E_0=\Fqt$ and permutes the $(\ell-1)$-tuple $E_1,\ldots,E_{\ell-1}$, since $x\in E_L$ implies $x^p\in E_{Lp}$. Let $E_1\ii=E_M$, hence $\eta\ii\in E_M$, $M\neq0$.
Since \eqref{e:sonoequiv2} also holds for $u_1=u_2=\ldots=u_{h-1}=u_{h+1}=\ldots=u_k=0$, for any $u=u_h\in I$ there exists $x=x_j\in\tilde I$ satisfying
\begin{equation}\label{e:v01_2}
\begin{cases}x&=\alpha_0u\ii+\alpha_{-M}\eta\ii{u\ii}^{q^s}\\ \tilde\eta x^{q^{\tilde s}}&=\alpha_1u\ii+\alpha_{1-M}\eta\ii{u\ii}^{q^s}\\
0&=\alpha_2u\ii+\alpha_{2-M}\eta\ii{u\ii}^{q^s}\\
\vdots&=\qquad\vdots\\
0&=\alpha_{\ell-1}u\ii+\alpha_{\ell-1-M}\eta\ii{u\ii}^{q^s},
\end{cases}
\end{equation}
where the indices are taken $\text{mod }\ell$.
Since \eqref{e:v01_2} must hold for all $u\in I$, the last $\ell-2$ equations give $\alpha_2=\ldots=\alpha_{\ell-1}=0$.
We now need to distinguish three cases: $M\notin\{1,\ell-1\}$, $M=1$, and $M=\ell-1$, respectively.

Case 1: $M\notin\{1,\ell-1\}$.\\
From \eqref{e:v01_2}, since $\alpha_{-M}=\alpha_{1-M}=0$,
\begin{equation}\label{e:24}
\begin{cases}x&=\alpha_0u\ii\\ \tilde\eta x^{q^{\tilde s}}&=\alpha_1u\ii.
\end{cases}
\end{equation}
Thus,
\[\tilde\eta\alpha_0^{q^{\tilde s}}u^{\iota q^{\tilde s}}-\alpha_1 u\ii=0.\]
Taking into account the assumption $\dim_\Fq I>2$ this implies $a_{jh}=0$ for all $j$ and $h$. Hence, this case cannot occur.

Case 2: $M=1$.\\ 
From \eqref{e:v01_2}, 
\begin{equation}
\begin{cases}x&=\alpha_0u\ii\\ \tilde\eta x^{q^{\tilde s}}&=\alpha_1u\ii+\alpha_0 \eta\ii u^{\iota q^{ s}}.
\end{cases}
\end{equation}
Thus,
\begin{equation}\label{e:v03_2}
\alpha_0^{q^{\tilde s}}\tilde\eta{u\ii}^{q^{\tilde s}}- \alpha_0\eta\ii{u\ii}^{q^s}-\alpha_1u\ii=0,
\end{equation}
for all $u\in I$.

Case 3: $M=\ell-1$.\\ 
From \eqref{e:v01_2}, 
\begin{equation}\label{eq:27}
\begin{cases}x&=\alpha_0u\ii+\alpha_1\eta\ii u^{\iota q^{ s}}\\ \tilde\eta x^{q^{\tilde s}}&=\alpha_1u\ii.
\end{cases}
\end{equation}
Thus,
\begin{equation}\label{e:v03_3}
\alpha_1^{q^{\tilde s}}\tilde\eta\eta^{\iota q^{\tilde s}}u^{\iota q^{s+\tilde s}}+\tilde\eta\alpha_0^{q^{\tilde s}}{u\ii}^{q^{\tilde s}}-\alpha_1u\ii=0,
\end{equation}
for all $u\in I$.

So far we have proved:\\
\emph{Assertion $(\star)$. Assume \eqref{e:sonoequiv2} is true, together with the assumption just above it. Then $\alpha_2=\alpha_3=\ldots=\alpha_{\ell-1}=0$, and either $(a)$~$E_1\ii=E_1$ and \eqref{e:v03_2}, or $(b)$~$E_1\ii=E_{\ell-1}$ and \eqref{e:v03_3}, hold for all $u\in I$.}

If $s\neq\tilde s,t-\tilde s$ and $I=\Fqt$, since both \eqref{e:v03_2} and \eqref{e:v03_3} admit less than $q^t$ solutions for $u$ unless all coefficients are zero, we obtain $a_{jh}=0$ for all $j$ and $h$, a contradiction. So, assertion $(i)$ is proved.

Now we prove $(ii)$.
Assume that $\left(T_{s,\eta,I}\right)^k$ and $\left(T_{s,\tilde\eta,\tilde I}\right)^k$ are equivalent under $\operatorname{\Gamma L}(k,q^{\ell t})$. 
Since in this case \eqref{e:v03_3} cannot be satisfied by all $u\in I$ unless all coefficients are zero, we are in case $(a)$ of Assertion $(\star)$.
Since  $s=\tilde s$, by \eqref{e:v03_2}, for all $u\in I$ the following equation holds:
\[
(\alpha_0^{q^{s}}\tilde\eta- \alpha_0\eta\ii){u\ii}^{q^s}-\alpha_1u\ii=0.
\]
Hence, $\alpha_1=0$ and $\alpha_0\tilde\eta(\alpha_0^{q^s-1}-\eta\ii\tilde\eta^{-1})=0$. It follows that for any $j,h\in\{1,\ldots,k\}$, $a_{jh}\in\Fqt$ and $a_{jh}(a_{jh}^{q^s-1}-\eta\ii\tilde\eta^{-1})=0$.
Since $A$ has some non-zero entries, this implies that $\N_{q^t/q}(\eta\ii\tilde\eta^{-1})=1$.
Furthermore, $a_{jh}\in\la \Theta\ra_{\Fq}$,  where $\Theta^{q^s-1}=\eta^\iota\tilde\eta^{-1}$, for all $j$ and $h$.
By \eqref{e:24},
$ a_{jh}I\ii\subseteq \tilde I$
and the equality must hold for the non-zero values among the $a_{jh}$s.

Conversely, assume that the conditions in $(ii)$ hold.
The equation
\[
\Theta \begin{pmatrix}u_1\ii+\eta\ii {u_1\ii}^{q^s}\\ u_2\ii+\eta\ii {u_2\ii}^{q^s}\\ \vdots\\u_k\ii+\eta\ii {u_k\ii}^{q^s}\end{pmatrix}=\begin{pmatrix}x_1+\tilde\eta x_1^{q^s}\\ x_2+\tilde\eta x_2^{q^s}\\ \vdots\\ x_k+\tilde\eta x_k^{q^s}\end{pmatrix},
\]
is equivalent to
\[
\begin{cases}
    x_j&=\Theta u_j\ii\\ \tilde\eta\Theta^{q^s}{u_j\ii}^{q^s}&=\eta\ii \Theta{u_j\ii}^{q^s},
\end{cases}\quad j=1,\ldots,k.
\]
Since the second equation is an identity, the semilinear map \[(\delta_1,\delta_2,\ldots,\delta_k)\mapsto \Theta(\delta_1\ii,\delta_2\ii,\ldots,\delta_k\ii),\] transforms $\left(T_{s,\eta,I}\right)^k$ into $\left(T_{s,\tilde\eta,\tilde I}\right)^k$.

The proof of assertion $(iii)$ is similar. Assume that $\left(T_{s,\eta,I}\right)^k$ and $\left(T_{t-s,\tilde\eta,\tilde I}\right)^k$ are equivalent under $\operatorname{\Gamma L}(k,q^{\ell t})$. We are now in case~$(b)$ of Assertion~$(\star)$. 
By Equation \eqref{e:v03_3}, for all $u\in I$ the following equation holds:
\[
(\alpha_1^{q^{t-s}}\tilde\eta\eta^{\iota q^{t-s}}-\alpha_1)u\ii+ \alpha_0^{q^{t-s}}\tilde\eta u^{\iota q^{t-s}}=0.
\]
This implies $a_{jh}\in E_1$, and $a_{jh}\left(a_{jh}^{q^{t-s}-1}\tilde\eta{\eta\ii}^{q^{t-s}}-1\right)=0$.
Since $\eta,\tilde\eta\in E_1$, we have $\eta\ii\in E_{\ell-1}$, ${\eta\ii}^{q^{t-s}}\in E_{-q^{t-s}}$ and $\tilde\eta{\eta\ii}^{q^{t-s}}\in E_{1-q^{t-s}}$.
Furthermore, there exists $\Lambda\in E_1$ such that $a_{jh}\in\la\Lambda\ra_{\mathbb F_q}$ for all $j,h=1,\ldots,k$ and $\Lambda^{1-q^{t-s}}=\tilde\eta\eta^{\iota q^{t-s}}$.
From the first equation in System~\eqref{eq:27} we have $\tilde I=\Lambda\eta\ii {I\ii}^{q^s}$.

The converse in~$(iii)$ can be proved in a way similar to the previous part.    
\end{proof}


\section{The regular fat polynomial \texorpdfstring{$\phi_{m,\sigma}$}{phi	extsubscript{m,sigma}}}\label{s:4}

In this section, the properties of the $\Rt$partially scattered polynomial $\phi=\phi_{m,\sigma}$ in \eqref{q:pol} are investigated.
The current assumptions are that $q$ is odd and that $\sigma=q^J$ with $\gcd(J,2t)=1$.
Let \[E=\{x\in\F_{q^{2t}}\colon x^{q^t}+x=0\},\] 
that is $E=E_1$ as defined in \eqref{e:ej} for $\ell=2$, and $E^*=E\setminus\{0\}$.
The main theorem in \cite{SmZaZu2} leaves open the question of the weights of the points in $L_\phi$ when $m$ is either a $(\sigma+1)$-th or a $(\sigma-1)$-th power of an element of $E$ (\footnote{In fact, $(\sigma-1)$-powers of elements in $E$ are $(q-1)$-powers and $(\sigma+1)$-powers are $(q+1)$-powers, and conversely.}).
In fact, these are two separate cases.
\begin{thm}\label{powers}
  The set of $(\sigma+1)$-th powers and the set of $(\sigma-1)$-th powers of elements in $E^*$ are disjoint.
\end{thm}
\begin{proof}
    Fix a generator $\omega$ of the multiplicative group $(\F_{q^{2t}}^*,\cdot)$, then \begin{equation}\label{e:struttE}E^*=\left\{\omega^{\frac{(1+2i)(q^t+1)}{2}}\colon i=1,\ldots,q^t-1 \right\}.\end{equation} 
    Take by contradiction $y=\omega^{\frac{(1+2i)(q^t+1)}{2}},\, z=\omega^{\frac{(1+2j)(q^t+1)}{2}}\in E^*$ such that $y^{\sigma+1}=z^{\sigma-1}$. Hence,
    \[q^{2t}-1 \mid  \frac{(\sigma+1)(1+2i)(q^t+1)}{2}-\frac{(\sigma-1)(1+2j)(q^t+1)}{2},\]
    which implies
    \begin{equation}\label{divisors}
      q^t-1\mid1+(\sigma+1)i-(\sigma-1)j.  
    \end{equation}
    We get a contradiction since in Relation~\eqref{divisors} the term on the left is even and the term on the right is odd, for any choice of $i$ and $j$.
\end{proof}

\subsection{Case 1: \texorpdfstring{$m$ is a $(\sigma+1)$-th power of an element of $E$}{m is a (sigma+1)-th power of an element of E}}

In this section we show that in this case $L_\phi$ is a linear set with complementary weights.

\begin{prop}\label{p:npsz}
\cite[Corollary 4.9]{NapPolSanZul22}
Let $s$ be a positive integer coprime with $t$, $\mu\in\F_{q^t}^*$ and $\xi\in\F_{q^{2t}}\setminus\Fqt$ such that $\N_{q^t/q}(\mu)\neq1$ and $\N_{q^t/q}(-\xi^{q^t+1}\mu)\neq(-1)^t.$
Let
\begin{equation}\label{e:nspz}
U=\{(u+\xi u^{q^s},v+\xi\mu v^{q^s})\colon u,v\in\Fqt\},
\end{equation}
then $L_U$ has exactly two points of weight greater than one.
\end{prop}


\begin{thm}\label{t:s+1}
Let $m$ be a $(\sigma+1)$-power of an element of $E^*$.
\begin{enumerate}[$(i)$]
\item If $t$ is odd, then $L_\phi$ is projectively equivalent to one of the linear sets described in Proposition~\ref{p:npsz}. Hence, it is a $(2,t)$-regular fat linear set in $\PG(1,q^{2t})$.
\item If $t$ is even, then $L_\phi$ is projectively equivalent to one of the linear sets described in Theorem~\ref{t:mainexample1}. Hence, it is a $(q+1,t)$-regular fat linear set in $\PG(1,q^{2t})$.
\end{enumerate}
\end{thm}
\begin{proof}
Since $E^*$ is closed with respect to inversion, there exists $w\in E^*$ such that $m=w^{-(\sigma+1)}$.
Let $x=x_0+wx_0^{\sigma^{t-1}}$, $y=y_0-wy_0^{\sigma^{t-1}}$, $x_0,y_0\in\Fqt$, and \[A=\begin{pmatrix}2&w\\ 2&-w\end{pmatrix}.\]
Then, using
\[
\phi(z)=\begin{cases}2z^{\sigma^{t-1}}&\text{ for }z\in\Fqt,\\
2w^{-(\sigma+1)}z^\sigma&\text{ for }z\in E,\end{cases},
\]
\begin{eqnarray*}
A\begin{pmatrix}x\\ \phi(x)\end{pmatrix}&=\begin{pmatrix}(2x_0+2wx_0^{\sigma^{t-1}})+w(2x_0^{\sigma^{t-1}}+2w^{-(\sigma+1)}w^\sigma x_0)\\
(2x_0+2wx_0^{\sigma^{t-1}})-w(2x_0^{\sigma^{t-1}}+2w^{-(\sigma+1)}w^\sigma x_0)\end{pmatrix}=
\begin{pmatrix}4x\\0\end{pmatrix},\\
A\begin{pmatrix}y\\ \phi(y)\end{pmatrix}&=\begin{pmatrix}(2y_0-2wy_0^{\sigma^{t-1}})+w(2y_0^{\sigma^{t-1}}-2w^{-(\sigma+1)}w^\sigma y_0)\\
(2y_0-2wy_0^{\sigma^{t-1}})-w(2y_0^{\sigma^{t-1}}-2w^{-(\sigma+1)}w^\sigma y_0)\end{pmatrix}=\begin{pmatrix}0\\ 4y\end{pmatrix}.
\end{eqnarray*}
Hence $AU_\phi$ contains the $t$-dimensional $\Fq$-subspaces
\begin{eqnarray*}
T_1&=\left\{(x_0+wx_0^{\sigma^{t-1}},0)\colon x_0\in\Fqt\right\},\\
T_2&=\left\{(0,y_0-wy_0^{\sigma^{t-1}})\colon y_0\in\Fqt\right\}.
\end{eqnarray*}
This implies $AU_\phi=T_1\oplus T_2$.
So, $U_{\phi}$ is in the same orbit of 
\begin{equation}\label{e:fcw}
T_{J(t-1),w,\Fqt}\ \times\ \left\{y_0-wy_0^{\sigma^{t-1}}\colon y_0\in\Fqt\right\},
\end{equation}
(cf.\ \eqref{e:defT}) under the action of $\GL(2,q^{2t})$.

If $t$ is odd, then $s=J(t-1)$, $\xi=w$ and $\mu=-1$ satisfy the assumptions of Proposition~\ref{p:npsz}. 
In particular
\[
\N_{q^t/q}(-\xi^{q^t+1}\mu)=\N_{q^t/q}(-w^2)\neq-1,
\]
by Proposition~\ref{p:haresatisfied}.
So, the $\Fq$-subspaces in Expressions~\eqref{e:nspz} and \eqref{e:fcw} coincide.

If $t$ is even, then $\N_{q^t/q}(-1)=1$ and a $\nu\in\Fqt$ exists such that $\nu^{\sigma^{t-1}-1}=-1$.
In this case, the substitution $y_0=\nu z$, $z\in\Fqt$ leads to
\[
y_0-wy_0^{\sigma^{t-1}}=\nu(z+wz^{\sigma^{t-1}}),
\]
and the statement in $(ii)$ is proved.
\end{proof}

\subsection{Case 2: \texorpdfstring{$m$ is a ($\sigma$-1)-th power of an element of $E$}{m is a (sigma-1)-th power of an element of E}}

In this subsection, we assume $m=w^{\sigma-1}$ for $w\in E^*$.
 The weight of a point $P_x=(x:\phi(x))$, $x\in\F_{q^{2t}}^*$, is $d$ if and only if the equation $\phi(\rho x)=\rho\phi(x)$ in the unknown $\rho$ has precisely $q^d$ solutions. Let $x=x_0+x_1$, $\rho=h+r$, $x_0,h\in \Fqt$, $x_1,r\in E$.
 That equation is equivalent to the following system in the unknowns $h\in\Fqt$ and $r\in E$ \cite{SmZaZu2}:
\begin{equation}\label{eq:systrrq}
\left\{
\begin{array}{lll}
-x_1r+m^{\sigma}x_1^{{\sigma}^2} r^{\sigma}&=&(h-h^{\sigma})x_0,\\
-x_0^{{\sigma}^{t-1}}r+mx_0^{\sigma}r^{\sigma}&=&m(h-h^{\sigma})x_1^{\sigma}.
\end{array}
\right.
\end{equation}

Our next goal is to find properties of the points of weight greater than one.
\begin{prop}\label{p:pc0}
    Any point $P_x$ with $x\in\F_{q^t}^*$ or $x\in E^*$ has weight two.
\end{prop}
\begin{proof}
Assume $x\in\F_{q^t}^*$, so $x_0=x$, $x_1=0$.
Equations~\eqref{eq:systrrq} imply $h\in\Fq$, and
\[
-x_0^{{\sigma}^{t-1}}r+mx_0^{\sigma}r^{\sigma}=0,
\]
that is, $r=0$ or $r^{\sigma-1}=w^{1-\sigma}x_0^{\sigma^{t-1}-\sigma}$. Hence, System~\eqref{eq:systrrq} has $q^2$ solutions
\[
(h,r)=(h,\lambda w^{-1}x_0^{({\sigma^{t-1}-\sigma})/(\sigma-1)}),\quad h,\,\lambda\in\Fq.
\]

Let $x\in E^*$. Then $h\in\mathbb F_q$ as above. The first equation in~\eqref{eq:systrrq} reads
    \[-x_1r+w^{\sigma^2-\sigma} x_1^{\sigma^2}r^\sigma=0,\]
    i.e., $r=0$ or $r^{\sigma-1}=w^{\sigma-\sigma^2}x_1^{1-\sigma^2}=(w^{-\sigma}x_1^{-\sigma-1})^{\sigma-1}.$ This means that the ratio between $r$ and $w^{-\sigma}x_1^{-\sigma-1}$ is in $\fq$, leading to the same conclusion.
\end{proof}
\begin{rem}
The previous proposition involves two sets of weight-two points contained in $L_\phi$,
\[
\{(x:\phi(x))\colon x\in\F_{q^t}^*\},\quad
\{(x:\phi(x))\colon x\in E^*\}.
\]
In fact, they are both of size $(q^t-1)/(q-1)$, and coincide.
Indeed, the map $x\mapsto w^{-1}x^{-\sigma^{t-1}}$ is a bijection between $\F_{q^t}^*$ and $E^*$, and for $x\in\F_{q^t}^*$
\[
\frac{\phi(w^{-1}x^{-\sigma^{t-1}})}{w^{-1}x^{-\sigma^{t-1}}}=2x^{\sigma^{t-1}-1}=\frac{\phi(x)}x.
\]
\end{rem}

System~\eqref{eq:systrrq} can be seen as a linear system in the unknowns $r$ and $r^\sigma$, and the determinant of the coefficient matrix is
\begin{equation}\label{e:D}
D=\begin{vmatrix}-x_1 & m^{\sigma}x_1^{{\sigma}^2} \\
-x_0^{{\sigma}^{t-1}} & m x_0^{\sigma}\end{vmatrix}.
\end{equation}

\begin{prop}\label{p:pc1}
    Consider a point $P_x=(x:\phi(x))$, where $0\neq x=x_0+x_1$, $x_0\in\F_{q^t}$, $x_1\in E$, and $D=0$, then $w_{L_\phi}(P_x)=2$.
\end{prop}
\begin{proof}
If $x_0x_1=0$, then the point has weight two by Proposition~\ref{p:pc0}. Then, assume $x_0x_1\neq0$. 
If $h-h^\sigma\neq0$, as shown in the proof of \cite[Theorem 2.3]{SmZaZu2}, then, because of the consistency of System~\eqref{eq:systrrq}, also
\[
\begin{vmatrix}-x_1&x_0\\ -x_0^{\sigma^{t-1}}&mx_1^\sigma\end{vmatrix}=0,
\]
and
$m=\left(x_0^{\sigma^{t-1}}x_1^{-1}\right)^{\sigma+1}$ would be a $(\sigma+1)$-th power of an element of $E^*$, contradicting Theorem \ref{powers}. Hence, $h\in\fq$. 
Since the determinant of the coefficient matrix is 0, the system \eqref{eq:systrrq} is equivalent to:
\begin{equation}\label{eq:systrrq_1}
   h\in\Fq,\quad-x_1 r+m^\sigma x_1^{\sigma^2}r^\sigma=0.
\end{equation} We obtain $r=0$, or
    \[r^{\sigma-1}=m^{-\sigma}x_1^{1-\sigma^2}=(w^{-\sigma}x_1^{-\sigma-1})^{\sigma-1},\]
    so $r=\xi w^{-\sigma}x_1^{-\sigma-1}$, $\xi\in\fq$. Thus, the number of solutions $(h,r)$ is $q^2$.
\end{proof}

Now we complete the proof that $L_\phi$ is a regular fat linear set.

\begin{thm}\label{t:s-1}
Every point of $L_\phi$ has weight at most two.
\end{thm}
\begin{proof}
By Propositions~\ref{p:pc0} and \ref{p:pc1}, $Dx_0x_1\neq0$ may be assumed.
From System~\eqref{eq:systrrq},
\[
mx_1^\sigma(-x_1r+m^{\sigma}x_1^{{\sigma}^2} r^{\sigma})-x_0(-x_0^{{\sigma}^{t-1}}r+mx_0^{\sigma}r^{\sigma})=0
\]
which is an equation in $r$, $r^\sigma$ that does not vanish because of $D\neq0$.
By Theorem~\ref{t:GQ}, there are at most $q$ possible solutions for $r$ and they belong to a one-dimensional $\Fq$-subspace of $E$.
For fixed $r$, each of the equations in~\eqref{eq:systrrq} is of type $d(h-h^\sigma)+d'=0$ with $d\neq0$, having at most $q$ solutions for $h$.
\end{proof}

We provide a characterization of the points of weight two for $L_\phi$. Since this is not essential for our purposes, we will omit the proof.

\begin{thm}
    Let $m=w^{\sigma-1}$ and $0\neq x=x_0+x_1$, where $x_0\in \Fqt$, $w,x_1\in E$.
    Then the point $P_x$ has weight two for $L_{\phi}$ if and only if 
    \[
    \Tr_{q^t/q}\left(\frac {-mx_0^\sigma x_1+m^\sigma x_0^{\sigma^{t-1}}x_1^{\sigma^2}}{w^\sigma (x_0^{\sigma^{t-1}+1}-mx_1^{\sigma+1})^{\sigma+1}}\right)=0.
    \]
\end{thm}

\subsection{Summary}
The following result summarizes \cite[Theorem 2.3]{SmZaZu2}, Theorem~\ref{t:s+1}, and Theorem~\ref{t:s-1}, and states that $L_{\phi_{m,\sigma}}$ is a regular fat linear set for all $m\in\F_{q^t}^*$.

\begin{thm}\label{t:summary}
Let
$
{\phi_{m,\sigma}}=
X^{\sigma^{t-1}}+X^{\sigma^{2t-1}}+m\left(X^\sigma-X^{\sigma^{t+1}}\right)\in\F_{q^{2t}}[X]
$, $t\ge3$, $q$ odd, $\sigma=q^J$, $\gcd(J,2t)=1$, and $m\in\F_{q^t}^*$.
\begin{itemize}
    \item If $m$ is a $(\sigma-1)$-power of an element of 
    {$E=\ker\operatorname{Tr}_{q^{2t}/q^t}$,}
    then $L_{\phi_{m,\sigma}}$ is an $(r,2)$-regular fat linear set for some integer $r$. 
    \item If $m$ is a $(\sigma+1)$-power of an element of $E$ and $t$ is odd, then $L_{\phi_{m,\sigma}}$ is a $(2,t)$-regular fat linear set.
    \item If $m$ is a $(\sigma+1)$-power of an element of $E$ and $t$ is even, then $L_{\phi_{m,\sigma}}$ is a $(q+1,t)$-regular fat linear set.
        \item Otherwise, $L_{\phi_{m,\sigma}}$ is a  scattered linear set.
\end{itemize}
\end{thm}

\section{Codes associated with regular linear sets}\label{s:5}

We will present rank-metric codes in the ambient space $\mathbb{F}_{q^n}^m$, which is isometric to those of the $n\times m$ matrices over $\fq$.

The \emph{rank} (weight) $w(v)$ of a vector $v=(v_1,v_2,\ldots,v_m) \in \mathbb{F}_{q^n}^m$ is defined as $w(v)=\dim_{\fq} (\langle v_1,v_2,\ldots, v_m\rangle_{\fq})$, {which equals the rank of the $n\times m$ matrix of the coordinates of $v_1,v_2,\ldots, v_m$ with respect to an $\fq$-basis of $\mathbb F_{q^n}$.}

A \emph{(linear vector) rank-metric $[m,k,d]_{q^n/q}$ code} $\C $ is a $k$-dimensional $\F_{q^n}$-subspace of $\mathbb{F}_{q^n}^m$ such that \[ d=\min\{w(v)\colon v\in\C,\,v\neq0\}.\] It is endowed with the rank distance defined as $d(x,y)=w(x-y)$ for all $x, y \in \mathbb{F}_{q^n}^m$.

\begin{thm}(Rank-metric Singleton bound)\label{th:singletonrank}
    Let $\C$ be an $[m,k,d]_{q^n/q}$ code.
Then 
\begin{equation}\label{eq:boundgen}
nk \leq \max\{m,n\}(\min\{m,n\}-d+1).\end{equation}
\end{thm}

An $[m,k,d]_{q^n/q}$ code $\C$ is said to be \emph{non-degenerate} if the $\fq$-span of the columns of any $k\times m$ generator matrix $G$ of $\C$ has $\fq$-dimension $m$.
This is equivalent to the following.

\begin{prop}\cite[Proposition 3.2]{lincutt}\label{prop:dCperpatleast2}
    Consider an $[m,k,d]_{q^n/q}$ code $\C$. The code $\C$ is non-degenerate if and only if the minimum distance of $\C^\perp$ is at least two, where $\C^\perp$ is the dual of $\C$ with respect to the standard dot product.
\end{prop}

We can now recall the connection between rank-metric codes and systems. 

\begin{thm}(see \cite{Randrianarisoa2020ageometric}) \label{th:connection}
Let $\C$ be a non-degenerate $[m,k,d]_{q^n/q}$ rank-metric code and let $G$ be a generator matrix of $\C$.
Let $U \subseteq \F_{q^n}^k$ be the $\F_q$-span of the columns of $G$.
The rank of an element $x G \in \C$, with $x \in \F_{q^n}^k$ is
\begin{equation}\label{eq:relweight}
w(x G) = m - \dim_{\fq}(U \cap x^{\perp}),\end{equation}
where $x^{\perp}=\{y \in \F_{q^n}^k \colon x \cdot y=0\}.$ In particular,
\begin{equation} \label{eq:distancedesign}
d=m - \max\left\{ \dim_{\fq}(U \cap H)  \colon H\mbox{ is an } \F_{q^n}\mbox{-hyperplane of }\F_{q^n}^k  \right\}.
\end{equation}
\end{thm}

Any non-degenerate code can be studied through an associated system. 
An $[m,k,d]_{q^n/q}$ \emph{system} $U$ is an $\mathbb{F}_q$-subspace of $\mathbb{F}_{q^n}^k$ of dimension $m$ such that 
\[
\langle U \rangle_{\mathbb{F}_{q^n}} = \mathbb{F}_{q^n}^k,
\]
and
\[
d = m - \max \left\{ \dim_{\mathbb{F}_q}(U \cap H) \colon H \text{ is an } \mathbb{F}_{q^n}\text{-hyperplane of } \mathbb{F}_{q^n}^k \right\}.
\]

The above result allows us to establish a one-to-one correspondence between equivalence classes of non-degenerate 
$[m,k,d]_{q^n/q}$ codes and the equivalence classes (under the action of $\mathrm{GL}(k,q^n)$) 
of $[m,k,d]_{q^n/q}$ systems; see~\cite{Randrianarisoa2020ageometric} for details. 
The system $U$ and the code $\mathcal{C}$ appearing in Theorem~\ref{th:connection} are said to be \emph{associated}.

Before we use the above connection, we recall the following notion of duality.

Let $\mathbb V=\F_{q^n}^k$. Let $\sigma \colon \mathbb{V}\times \mathbb{V} \rightarrow \mathbb{F}_{q^n}$ be a non-degenerate reflexive sesquilinear form over $\mathbb{V}$. Define
$\sigma' \colon \mathbb{V} \times \mathbb{V} \rightarrow \mathbb{F}_q$ by $\sigma':(u,v)\mapsto {\Tr}_{q^n/q}(\sigma(u,v))$.
If we regard $\mathbb{V}$ as an $\F_q$-vector space, then $\sigma^\prime$ turns out to be a non-degenerate reflexive sesquilinear form on $\mathbb{V}$.
Let $\perp$ and $\perp'$ be the orthogonal complement maps defined by $\sigma$ and $\sigma'$ on the lattices of $\F_{q^n}$-linear and 
$\F_q$-linear subspaces, respectively.
The following properties hold (see \cite[Section~2]{polverino2010linear} for more details).

\begin{prop}\label{prop:dualityproperties}
With the above notation,
\begin{itemize}
    \item[(i)] $\dim_{\F_{q^n}}(W)+\dim_{\F_{q^n}}(W^\perp)=k$, for every $\F_{q^n}$-subspace $W$ of $\mathbb{V}$.
    \item[(ii)] $\dim_{\F_{q}}(U)+\dim_{\F_{q}}(U^{\perp'})=nk$, for every $\F_{q}$-subspace $U$ of $\mathbb{V}$.
    \item[(iii)] For all $\F_q$-subspaces $T_1$, $T_2$ of $\mathbb{V}$, $T_1\subseteq T_2$ implies $T_1^{\perp'}\supseteq T_2^{\perp'}$.
    \item[(iv)] $W^\perp=W^{\perp'}$, for every $\F_{q^n}$-subspace $W$ of $\mathbb{V}$. 
\end{itemize}
\end{prop}

\begin{prop}\label{prop:dualconsidercode}
    Let $U$ be an $\fq$-subspace of $\F_{q^n}^k$ such that $L_U$ is an $(r,i)$-regular fat linear set in $\PG(k-1,q^n)$ with $i<n$. We have $\langle U^{\perp'}\rangle_{\F_{q^n}}=\F_{q^n}^k$.
\end{prop}
\begin{proof}
    Suppose that a hyperplane $H$ of $\mathbb V$ contains $U^{\perp'}$. Then, $H^\perp \subseteq U$. Since $\dim_{\F_q}(H^\perp)=n$, there is a point of weight $n$ in $L_U$, which is a contradiction.
\end{proof}

Therefore, the above proposition allows us to consider the code associated with the dual of an $(r,i)$-regular fat linear set. Note that the size of an $(r,i)$-regular fat linear set is given by~\eqref{e:size}.

\begin{prop}\label{prop:weightdistribution}
     Let $U$ be an $\fq$-subspace of {$\F_{q^n}^k$ and $\dim_{\Fq}(U)=\rho$, } such that $L_U$ is an $(r,i)$-regular fat linear set with $i<n$.
    Consider the rank-metric code $\C$ associated with $U^{\perp'}$.
   Then $\C$ is a rank-metric $[nk-\rho,k,n-i]_{q^n/q}$-code. Such $\C$ has
    \begin{itemize}
        \item $r(q^n-1)$ codewords of weight $n-i$;
        \item $(|L_U|-r)(q^n-1)$ codewords of weight $n-1$;
        \item $q^{nk}-1-|L_U|(q^n-1)$ codewords of weight $n$.
    \end{itemize}
\end{prop}
\begin{proof}
    Note that for any $\F_{q^n}$-subspace $S=\mathrm{PG}(W,\mathbb{F}_{q^n})$ of $\mathrm{PG}(k-1,q^n)$ we have that
    \[ w_{L_U}(S)=\dim_{\fq}(W\cap U)= nk-\dim_{\fq}((W\cap U)^{\perp'})\]
    \[=nk- \dim_{\fq}(W^{\perp'}+U^{\perp'})=nk-(nk-n\dim_{\fqn}(W))-(nk-\rho)+\dim_{\fq}(U^{\perp'}\cap W^{\perp'})\]
    \[=n\dim_{\mathbb{F}_{q^n}}(W)+\rho-nk+w_{L_{U^{\perp'}}}(S^\perp). \]
    In particular, for any point $P \in \mathrm{PG}(k-1,q^n)$ we have that
    \[ w_{L_{U^{\perp'}}}(P^\perp)=nk-\rho-n+w_{L_U}(P), \]
    implying that for $L_{U^{\perp'}}$ there are
    \begin{itemize}
        \item $r$ hyperplanes of weight $nk-\rho-n+i$;
        \item $|L_U|-r$ hyperplanes of weight $nk-\rho-n+1$;
        \item the remaining hyperplanes have weight $nk-\rho-n$.
    \end{itemize}
    Using Theorem \ref{th:connection} we obtain the assertion.
\end{proof}

We obtain a first bound on the rank of an $(r,i)$-regular fat linear set.

\begin{prop}
    Let $U$ be an $\fq$-subspace of $\F_{q^n}^k$ with dimension $\rho$ such that $L_U$ is $(r,i)$-regular fat with $i<n$ and $\rho\leq nk-n$.
     Then
     \[ \rho\leq nki/(i+1). \]
\end{prop}
\begin{proof}
    Consider the rank-metric code $\C$ associated with $U^{\perp'}$.
    Note that $\C$ is an $[nk-\rho,k,n-i]_{q^n/q}$ code.
    Since $nk-\rho\geq n$, applying Theorem \ref{th:singletonrank} yields
    \[ nk\leq (nk-\rho)(n-(n-i)+1). \]
    This implies the thesis.
\end{proof}
We will now use MacWilliams-like identities to obtain more precise bounds on regular fat linear sets. Recall that the \textit{Gaussian binomial coefficient} $\begin{bmatrix}
         r \\
         h \\
        \end{bmatrix}_{q}$ is the number of $h$-dimensional $\Fq$-subspaces of the vector space $V=\mathbb F_{q}^r$.
    If $h=0$, we have $\begin{bmatrix}
         r \\
         h \\
        \end{bmatrix}_{q}=1$, otherwise if $1\leq h\leq r$,
   \begin{equation}
    \begin{bmatrix}
         r \\
         h \\
        \end{bmatrix}_{q}=\frac{(q^{r}-1)(q^{r-1}-1)\cdots(q^{r-h+1}-1)}{(q^{h}-1)(q^{h-1}-1)\cdots(q-1)}, 
       \end{equation} 
   while $\begin{bmatrix}
         r \\
         h \\
        \end{bmatrix}_{q}=0$ whenever $h>r$. 

\begin{thm}\cite[Theorem 31]{ravagnani2016rank}\label{thm:MWident}
    Let $\C$ be a rank-metric code in $\F_{q^n}^N$. Let $(A_i)_{i \in \mathbb{N}}$ and $(B_j)_{j \in \mathbb{N}}$ be the rank distributions of $\C$ and $\C^\perp$, respectively. For any integer $0\leq \nu \leq n$ we have
    \[ \sum_{j=0}^{n-\nu} A_j \begin{bmatrix}n-j\\
    \nu\\
    \end{bmatrix}_{q}=\frac{|\C|}{q^{N\nu}}\sum_{j=0}^\nu B_j\begin{bmatrix}n-j\\
    \nu-j\\
    \end{bmatrix}_{q}. \]
\end{thm}

\begin{thm}\label{th:mainBound}
    Let $U$ be an $\fq$-subspace of $\F_{q^n}^k$ such that $L_U$ is $(r,i)$-regular fat of rank $\rho$ with $i<n$.
    We have \[ r \geq \frac{(q^{2\rho-nk}-1)\begin{bmatrix}n\\2\\
    \end{bmatrix}_q}{(q^n-1)\begin{bmatrix}i\\2\\
    \end{bmatrix}_q}. \]
\end{thm}
\begin{proof}
    Consider the rank-metric code $\C$ associated with $U^{\perp'}$.
    Recall that $\C$ is an $[nk-\rho,k,n-i]_{q^n/q}$ code with rank distribution given in Proposition \ref{prop:weightdistribution}.
    Clearly $B_0=1$ and $B_1=0$ (by Proposition \ref{prop:dCperpatleast2}). Let us now determine $B_2$ via Theorem \ref{thm:MWident} with $\nu=2$.
    In this case, we have
    \[ \sum_{j=0}^{n-2}A_j \begin{bmatrix}n-j\\2\\
    \end{bmatrix}_{q}= \frac{q^{nk}}{q^{2(nk-\rho)}}\left(\begin{bmatrix}n\\2\\
    \end{bmatrix}_{q}+B_2\right), \]
    from which we derive 
    \[ B_2=q^{nk-2\rho}\left( \sum_{j=0}^{n-2}A_j \begin{bmatrix}n-j\\2\\
    \end{bmatrix}_{q} - q^{2\rho-nk}\begin{bmatrix}n\\2\\
    \end{bmatrix}_{q}\right)\]
    \[=q^{nk-2\rho}\left( (q^n-1)r\begin{bmatrix}i\\2\\
    \end{bmatrix}_q + \begin{bmatrix}n\\2\\
    \end{bmatrix}_q - q^{2\rho-nk}\begin{bmatrix}n\\2\\
    \end{bmatrix}_{q}\right).\]
    Since $B_2\geq 0$, we have the assertion.
\end{proof}

\section{Conclusions}\label{s:conclusions}

In this work, we introduced the notions of $(r, i)$-regular fat linear sets and $(r, i)$-regular fat linearized polynomials. Regular fat linear sets form a broad class that includes scattered linear sets, clubs, and other families that have been studied previously. We analyzed explicit constructions, such as those arising from the polynomial $\phi_{m,\sigma}$. Furthermore, we illustrated the relationship between these geometric structures and the theory of rank-metric codes, especially the study of codes with few weights. Our results show that regular fat linear sets generate new instances of three-weight rank-metric codes and that this connection can be exploited to derive structural bounds on the associated
linear sets via MacWilliams identities. The general framework of regular fat linear sets opens up several avenues for research. First, it offers a unifying perspective on several well-known classes. Second, it raises natural questions about existence, classification, and equivalence.
Before concluding with a description of some particularly compelling problems for future investigation, we present the following result, which admits generalizations to dimensions greater than one and suggests that finding examples of $(r, i)$-regular fat linear sets with arbitrary $r$ could be difficult.

\begin{prop}
    If $L_U$ is an $(r,i)$-regular fat linear set of rank $2i$ in $\PG(1,q^n)$, then $r\le2$ or $r=q^j+1$ where $j$ is a divisor of $n$.
\end{prop}
\begin{proof}
    Assume $r>2$.
    Since $\PGL(2,q^n)$ is $3$-transitive on the points of $\PG(1,q^n)$, it may be assumed that $(1:0)$, $(0:1)$, and $(1:1)$ are points of weight $i$. Define
    \begin{eqnarray*}
        U_1=U\cap(1:0),&\text{\ }&T_1=\{x\in\Fqn\colon(x,0)\in U_1\},\\
        U_2=U\cap(0:1),&\text{\ }&T_2=\{y\in\Fqn\colon(0,y)\in U_2\}.
    \end{eqnarray*}
    Both $T_1$ and $T_2$ are $i$-dimensional $\Fq$-subspaces of $\Fqn$. This in turn implies $U=U_1\oplus U_2=T_1\times T_2$.
    Since $U\cap(1:1)\subseteq T_1\times T_2$ contains $q^i$ elements of type $(x,x)$, we have $T_1=T_2$.
    Let $S=\{b\in\Fqn\colon bT_1\subseteq T_1\}$, which is a subfield of $\Fqn$.

    Let $P=(1:\alpha)$ be a point in $\PG(1,q^n)$.
    The weight of $P$ with respect to $L_U$ is $i$ if and only if $(x,\alpha x)\in U$ for all $x\in T_1$; this is equivalent to $\alpha\in S$.
    Hence, the number of points of weight $i$ is $1+|S|$.
\end{proof}

To our knowledge, there are currently no known examples of $(r,i)$-regular fat linear sets of rank $2i$ with $r>q+1$.
This gives rise to the following problem:
\begin{problem}
    Find examples of, or prove the non-existence of, $(q^j+1,i)$-regular fat linear sets of rank $2i$ with $j>1$ in $\mathrm{PG}(1,q^n)$.
    (\footnote{After the manuscript was accepted for publication, information was received by the authors that the open problems 1 and 2 had been solved.})
\end{problem}


Regarding regular fat linear sets in projective spaces, {in Section~\ref{s:examples}} we have shown examples in $\mathrm{PG}(k-1,q^n)$ where the number of points of weight greater than one is $(q^k-1)/(q-1)$, and these points form a subgeometry, $\mathrm{PG}(k-1,q)$, of $\mathrm{PG}(k-1,q^n)$. 
As already mentioned in Section \ref{s:regular}, an example of $(k,t)$-regular fat linear set in $\mathrm{PG}(k-1,q^{2t})$ has been described in \cite{Zul23}. 
This leads to the following.

\begin{problem}
Construct $(r,i)$-regular fat linear sets in $\mathrm{PG}(k-1,q^n)$ whose set of points of weight greater than one has one of the following properties: 
\begin{enumerate}[1)]
    \item either the set has a cardinality that is strictly between $k$ and $(q^k-1)/(q-1)$, or  
    \item the set does not form a subgeometry $\mathrm{PG}(k-1,q)$ of $\mathrm{PG}(k-1,q^n)$.  
\end{enumerate}
\end{problem}

\section*{Acknowledgments}

The research was partially supported by the Italian National Group for Algebraic and Geometric Structures and their Applications (GNSAGA - INdAM). The research of the first two authors was supported by the project SID \emph{Linear Algebra methods in Combinatorics and Industrial Engineering: Coding Theory, Graphs and their interest in Assembly Line Problems} CUP C33C25001110005 of Dipartimento di Tecnica e Gestione dei Sistemi Industriali of the Università degli Studi di Padova.

\end{document}